\documentclass[letterpaper,acmtkdd]{acmsmall}
\let\orgsetcounter\setcounter
\usepackage{calc}

\acmArticle{}
\runningfoot{}
\gdef\firstfoot{}
\makeatletter
\ps@appheadings
\def\endbottomstuff{
\vskip-13pt
\strut
\end@float}
\makeatother

\usepackage{amsmath}
\usepackage{amssymb}
\usepackage{booktabs}
\usepackage[Q=yes,bsdiv=yes,div=yes]{examplep}
\usepackage{graphicx}
\usepackage[svgnames,dvinames,x11names]{xcolor}
\usepackage{listings}
\usepackage{textcomp}
\usepackage{xspace}
\usepackage{subfigure}
\usepackage{tabularx}
\usepackage{textcase}
\usepackage{url}

\usepackage{tkz-graph}
\usetikzlibrary{arrows,shapes} 

\graphicspath{{./}{./figures/}}

\colorlet{TufteRed}{red!80!black}

\newcommand{\la}{\lambda}

\newcommand{\pa}{\partial}

\newcommand{\no}{\nonumber}

\newcommand{\mb}{\mathbb}
\newcommand{\mf}{\mathfrak}

\newcommand{\f}{\frac}

\newcommand{\vI}{\vec{1}}

\newcommand{\dsp}{\displaystyle}
\newcommand{\order}{\mathcal{O}}

\newcommand{\namp}{{NetAlignMP}\xspace}
\newcommand{\namppp}{{NetAlignMP\texttt{++}}\xspace}
\newcommand{\mppp}{{MP\texttt{++}}\xspace}

  \def\clap#1{\hbox to 0pt{\hss#1\hss}}%
  \def\mathclap{\mathpalette\mathclapinternal}%
  \def\mathllapinternal#1#2{%
    \llap{$\mathsurround=0pt#1{#2}$}}%
  \def\mathrlapinternal#1#2{%
    \rlap{$\mathsurround=0pt#1{#2}$}}%
  \def\mathclapinternal#1#2{%
    \clap{$\mathsurround=0pt#1{#2}$}}%

  \providecommand{\mat}[1]{\boldsymbol{#1}}%
  \providecommand{\normof}[2][]{\left\| #2 \right\|_{#1}}%
  \DeclareMathOperator{\diag}{diag}%

  \providecommand{\itn}[1]{^{(#1)}}%
  \providecommand{\eps}{\varepsilon}%
  \providecommand{\kron}{\otimes}
  \DeclareMathOperator{\tvec}{vec}
  \providecommand{\pmat}[1]{\begin{pmatrix} #1 \end{pmatrix}}
  \providecommand{\bmat}[1]{\begin{bmatrix} #1 \end{bmatrix}}

  \DeclareMathOperator*{\minimize}{minimize}
  \DeclareMathOperator*{\maximize}{maximize}

  \providecommand{\subjectto}{\ensuremath{\text{\upshape subject to}}}
  \providecommand{\MINof}[1][]{{\displaystyle \minimize_{#1}}}

  \providecommand{\MAXof}[1][]{{\displaystyle \maximize_{#1}}}
  
  \providecommand{\MAXone}[3]{\begin{array}{ll} \MAXof[#1] & #2 \\ \subjectto & #3 \end{array}}
  \providecommand{\MAXtwo}[4]{\begin{array}{ll} \MAXof[#1] & #2 \\ \subjectto  & #3 \\ & #4 \end{array}}
  \providecommand{\MAXthree}[5]{\begin{array}{ll} \MAXof[#1] & #2 \\ \subjectto  & #3 \\ & #4 \\ & #5 \end{array}}

  \DeclareMathOperator{\E}{E}
  \DeclareMathOperator{\p}{\mathbb{P}}
  
  \DeclareMathOperator{\Cov}{Cov}

  \DeclareMathOperator{\Std}{Std}

  \providecommand{\set}{\mathcal}
  \providecommand{\graph}{\mathcal}

  \providecommand{\eps}{\varepsilon}

  \providecommand{\eye}{\mat{I}}
  \providecommand{\mA}{\ensuremath{\mat{A}}}
  \providecommand{\mB}{\ensuremath{\mat{B}}}
  \providecommand{\mC}{\ensuremath{\mat{C}}}
  \providecommand{\mD}{\ensuremath{\mat{D}}}
  
  \providecommand{\mF}{\ensuremath{\mat{F}}}

  \providecommand{\mL}{\ensuremath{\mat{L}}}

  \providecommand{\mP}{\ensuremath{\mat{P}}}
  \providecommand{\mQ}{\ensuremath{\mat{Q}}}
  
  \providecommand{\mS}{\ensuremath{\mat{S}}}
  
  \providecommand{\mU}{\ensuremath{\mat{U}}}
  
  \providecommand{\mW}{\ensuremath{\mat{W}}}
  \providecommand{\mX}{\ensuremath{\mat{X}}}
  \providecommand{\mY}{\ensuremath{\mat{Y}}}
  \providecommand{\mZ}{\ensuremath{\mat{Z}}}

  \providecommand{\vd}{\ensuremath{\vec{d}}}
  \providecommand{\ve}{\ensuremath{\vec{e}}}

  \providecommand{\vp}{\ensuremath{\vec{p}}}

  \providecommand{\vv}{\ensuremath{\vec{v}}}
  \providecommand{\vw}{\ensuremath{\vec{w}}}
  \providecommand{\vx}{\ensuremath{\vec{x}}}
  \providecommand{\vy}{\ensuremath{\vec{y}}}
  \providecommand{\vz}{\ensuremath{\vec{z}}}

\renewcommand{\vec}{\boldsymbol}
\renewcommand{\graph}{}
\renewcommand{\set}{}

\DeclareMathOperator*{\bound}{bound}

\renewcommand{\Box}{
\mathrel{\mathchoice{\mbox{\fontsize{7}{7}$\square$}}
                    {\mbox{\fontsize{6}{6}$\square$}}
                    {\mbox{\fontsize{5}{5}$\square$}}
                    {\mbox{\fontsize{5}{5}$\square$}}}}

\newcommand{\citet}{\citeN}
\newcommand{\Section}{\section}
\newcommand{\SubSection}{\subsection}

\newcommand{\acro}[1]{\textsc{\MakeTextLowercase{#1}}}

\newenvironment{inlinealgorithm}[2]%
{%
\vspace{\abovedisplayskip}%
\hrule%
\vspace{\belowdisplayskip}%
\noindent \textit{Algorithm} #1\\%
\noindent {\sffamily \small INPUT} #2\par%
\footnotesize%
}%
{%
\vspace{\abovedisplayskip}
\hrule
\vspace{\belowdisplayskip}
}

\usepackage{epstopdf}

\title{Message Passing Algorithms for Sparse Network
Alignment}

\author{%
MOHSEN BAYATI \affil{Stanford University}
DAVID F. GLEICH \affil{Purdue University}
AMIN SABERI \affil{Stanford University}
YING WANG \affil{Google}
}

\begin{abstract}
Network alignment generalizes and unifies several approaches
for forming a matching or alignment between the vertices of two graphs.
We study a mathematical programming framework for network alignment problem and a sparse variation of it where only a small number of matches between
the vertices of the two graphs are possible.
We propose a new message passing algorithm that allows us to
compute, very efficiently, approximate solutions to the sparse network alignment problems with graph sizes as large as hundreds of thousands of vertices.
We also provide extensive simulations comparing our algorithms with two of the best solvers for network alignment problems on two synthetic matching problems, two bioinformatics problems, and three large ontology alignment problems including a multilingual problem with a known labeled alignment.
\end{abstract}



\begin{document}
{\let\setcounter\orgsetcounter
\begin{bottomstuff}
Author's address: Mohsen Bayati, Graduate School of Business, Stanford, CA \newline
David F. Gleich, Computer Science Department, Purdue University, West Lafayette, IN \newline
Amin Saberi, Management Science and Engineering Department, Stanford, CA \newline
Ying Wang, Google, Mountain View, CA
\end{bottomstuff}
}

\maketitle

\Section{Introduction} The focus of this paper is to find
approximate isomorphisms, or alignments, between large graphs. This
problem is motivated by applications in several areas including
biology, computer vision, and natural language processing. For
example, the study of protein interactions across different species
has made
network alignment a common topic in computational biology \cite{singh2007-matching-topology,singh2008-isorank-multi,%
klau2009-network-alignment,flannick2006-graemlin,flannick2008-graemlin2,kuchaiev2009-topological-alignment}.
In computer vision, network alignment is used for matching images
\cite{conte2004-graph-matching,scheelewald2005-subgraph-matching},
and in the ontology alignment, it is used for finding correspondence
between different representations of a database
\cite{melnik2002-similarity-flooding,svab2007-exploiting,Lacoste-Julien-2006-word-alignment}.

The formulation of the problem studied in this paper is a variation
of classic algorithmic problems: graph isomorphism,  maximum common
subgraph, and the quadratic assignment problem.
Because of the intractability of the problem, our  focus will be on
practical heuristics. We will give a quick review of the existing
results and their applications. Then we will present two message
passing algorithms that yield near optimal results -- determined
by comparison to an upper bound from a linear program.  Both
algorithms easily work on graphs with 100,000-1,000,000 vertices.
Because our algorithms use message passing, they
can be parallelized on MapReduce and bulk-synchronous
processing architectures for even larger problems.

\subsection{Problem Definition}

Consider two sets of vertices $V_A=\{1,2,\ldots,n\}$ and
$V_B=\{1',2',\ldots,m'\}$. Let $A = (V_A,E_A)$ and $B = (V_B, E_B)$
be two undirected graphs with their respective vertex and edge sets.
Let $L$ be a bipartite graph between the vertices of $A$ and $B$,
formally $L = (V_A \cup V_B, E_L)$.  Our overall goal is to find a
matching between $A$ and $B$ using only edges from $L$. In other
words, we seek a subset of $E_L$ such that no two edges share a
common endpoint. Under such a matching $M$, we say that an edge
$(i,j) \in E_A$ is \emph{overlapped} with $(i',j') \in E_B$ if
$(i,i')$ and $(j,j')$ belong to $M$.  See
Figure~\ref{fig:netalign-setup} for an illustration.

\begin{narrowfig}{0.5\columnwidth}[t]
   \includegraphics[width=0.5\columnwidth]{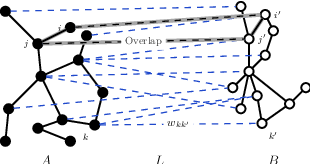}
   \caption{The setup for the network alignment problem.  The goal
   is to maximize the number of overlaps
   in any matching subset of $L$ and the weight of the matching.}
   \label{fig:netalign-setup}
 \end{narrowfig}

More generally, and following~\citet{singh2007-matching-topology}, we
will study the case where the edges between $A$ and $B$ are
weighted. That is, each edge $e=(k,k') \in E_L$ has a non-negative
weight $w_e$ indicating a measure of similarity between vertices $k$
and $k'$.  In these cases, a matching has a weight that is equal to
the sum of the weights of edges in the matching.

\begin{definition} Given graphs $\graph{A}$, $\graph{B}$ and
$\graph{L}$, as well as the weight function $w$, find a matching $M$
maximizing a linear combination of the matching weight and the
number of overlapped edges.
\end{definition}

The above problem is a generalization of several NP-complete
problems including the densest subgraph problem as well as the
maximum common subgraph problem. The latter is also known to be
APX-hard.

\subsection{Our contribution}

In this paper we provide
\begin{enumerate}
 \item two novel message passing algorithms: \namp and \namppp for the problem based on max-product belief propagation (Section \ref{sec:bp});
 \item an extensive comparison between \namp, \namppp and two of the best existing algorithms on two synthetic matching problems (Section \ref{sec:synthetic}), two bioinformatics problems, and three large ontology alignment problems (Section \ref{sec:realdata}) including a multilingual problem with a known alignment.
\end{enumerate}

We will show that our algorithms are fast, robust, and yield
near-optimal\footnote{For the families of graphs that we study, we have a theoretical upper bound provided by linear programming for the objective function. Hence we can check the quality of all algorithms' solutions.} objective values for a large family of graphs, including real datasets. In one of the cases where there is a
known alignment produced by experts, our algorithms
recover a large fraction of the correct matches quickly, without any
tweaking.

In evaluating our algorithms on synthetic datasets, we will observe
that for both sparse and dense cases, our algorithms produce
near-optimal solution (using the theoretical upper bounds). As the number
of edges in $L$ increases, our approaches outperform existing methods
by a factor of 2 or more. On the other hand, when there are only
a small number of potential matches in $L$,
our results nearly ties with Klau's algorithm~\cite{klau2009-network-alignment}.

All of our algorithms are implemented in \textsc{Matlab} and the
software and datasets for this paper are available to public from
the web-page:
\begin{center}
\url{http://www.cs.purdue.edu/homes/dgleich/codes/netalign}	
\end{center}
We also include all the experimental code to reproduce the figures in
this paper. Using our network alignment \textsc{Matlab} package,
solving the network alignment problem in
Figure~\ref{fig:netalign-setup} is done with the following code.
\begin{lstlisting}[language=Matlab,%
  basicstyle={\ttfamily\scriptsize},keywordstyle=\mdseries,%
  upquote=true,xleftmargin=5mm,mathescape=false,columns=fixed]
load('data/example_overlap.mat');   % load the data
[S,w,li,lj]=netalign_setup(A,B,L);  % setup the quadratic program
x=netalignmp(S,w,1,1,li,lj);        % call our netalignmp
[mi ma mb]=mwmround(x,li,lj);       % round to matching
\end{lstlisting}

This manuscript is an extension of our previous paper
\cite{bayati2009-network-alignment} and includes
a full derivation of our belief propagation algorithm,
as well as a new algorithm.  It also includes a more
thorough experimental evaluation.

\section{A mathematical program for network alignment}\label{sec:math-programming}

In this section, we formulate the network alignment as a quadratic program (QP).
Let us start by introducing some notation in Table~\ref{tab:notation}.

\begin{table}[h]
\tbl{Notation for the paper.}{
\begin{tabular}{rl}
\toprule
$\graph{A},\set{S}$ & capital letters are sets and graphs\\
$\mA,\mB,\mS$ & bold capitals are matrices\\
$\vx$, $\vw$ & lowercase bold letters are vectors\\
$A_{ij}, \mS[ii',jj']$ & subscripts or brackets denote matrix entries
\\ \midrule
$(i,i')\,, ii'$ & edges in $\graph{L}$  \\
$ii' \Box jj'$  & squares in $E_L\times E_L$ \\
$\vI_n$ & $n$ by 1 vector of all ones \\
$\mA \vx\,, \mA\mB$ & standard vector and matrix products  \\
$\mA \bullet \mB=\sum_{ij} \mA_{ij} \mB_{ij}$ & matrix inner product
\\ \bottomrule
\end{tabular}
}
\label{tab:notation}
\end{table}

Given $A = (V_A, E_A)$, $B = (V_B, E_B)$, and $L = (V_A \cup
V_B, E_L)$, our goal is to produce a matching $M$ to maximize a
linear combination of overlap and matching weight. For each edge of
$E_L$ we will use the notations $(i,i')$ and $ii'$ interchangeably.
Each matching in $E_L$ is represented by a zero-one vector by
assigning a variable $x_{ii'}$ to each $ii' \in E_L$, which is equal
to $1$ if $ii'$ is in the matching or $0$ if it is not. For
convenience of notation, we define an ordering $\order_L$ over the
set $E_L$. We will use the same ordering in vector representation of
the edges. Note that $x_{ii'}$ is \emph{only defined} for edges in $L$,
and $|E_L| \ll |V_A|\cdot|V_B|$.
This differs from many other formulations of the problem where the
set $L$ is implicitly the full bipartite collection.

Next, we define a zero-one matrix $\mS$ of size
$|E_L|\times |E_L|$ indexed by edges of $E_L$.  Denote the entry
at row $ii'$ and column $jj'$ by $\mS[ii',jj']$, where
\[
\mS[ii',jj'] = \begin{cases}
                1 & \text{if $(i,j) \in E_A$ and $(i',j') \in E_B$}\\
                0 & \text{otherwise}.
               \end{cases}
\]
We also say that two edges $ii'$ and $jj'$ in $E_L$ form a
\emph{square} if $\mS[ii',jj']=1$ and denote it by $ii'\Box jj'$. In
other words, $\mS$ is the indicator matrix of all squares.


Let $\vx$ be the indicator vector for a matching.
The total number of overlapped edges is
\[(1/2) \vx^T \mS \vx=\sum_{\mathclap{ii'\Box jj'}} x_{ii'} x_{jj'}.\]
Moreover, let $w_{ii'}$ be the weight of each $ii'$ and denote the vector of all weights by $\vw$.
The constraint that $\vx$ must be a valid matching can be written by set of linear inequalities.
For all vertices $(i,i')\in\graph{L}$,
\[
  \sum_{\mathclap{j' : (ij') \in E_L}} x_{ij'} \le 1, \qquad
  \sum_{\mathclap{j : (ji') \in E_L}} x_{ji'} \le 1 \qquad
  x_{ii'} \in \{0, 1\}.\]
To write these constraints more compactly, define $\mC$ to be the
binary incidence matrix of graph $\graph{L}$ of dimensions $|V_L|\times |E_L|$.
Then the matching constraints can be written as $\mC \vx \leq \vI_{|V_{L}|}$.

\begin{figure*}[t]
 \centering
 \parbox{3.05in}{\includegraphics{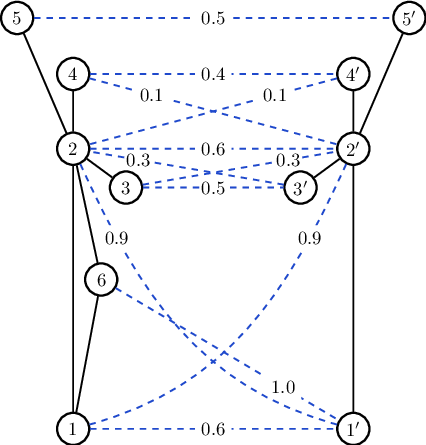}}

$  \begin{smallmatrix}
  \smash{\mathclap{\text{edge order $\order_L$}}} \\[0.5ex]
  (2,2')   \\
  (2,1')   \\
  (2,3')   \\
  (2,4')   \\
  (1,2')   \\
  (1,1')   \\
  (3,2')   \\
  (3,3')   \\
  (4,2')   \\
  (4,4')   \\
  (5,5')   \\
  (6,1')   \\
 \end{smallmatrix}
 \;\;
   \underbrace{\left(\begin{smallmatrix}
   0\vphantom{(0')} & 0 & 0 & 0 & 0 & 1 & 0 & 1 & 0 & 1 & 1 & 1 \\
   0\vphantom{(0')} & 0 & 0 & 0 & 1 & 0 & 1 & 0 & 1 & 0 & 0 & 0 \\
   0\vphantom{(0')} & 0 & 0 & 0 & 1 & 0 & 1 & 0 & 1 & 0 & 0 & 0 \\
   0\vphantom{(0')} & 0 & 0 & 0 & 1 & 0 & 1 & 0 & 1 & 0 & 0 & 0 \\
   0\vphantom{(0')} & 1 & 1 & 1 & 0 & 0 & 0 & 0 & 0 & 0 & 0 & 1 \\
   1\vphantom{(0')} & 0 & 0 & 0 & 0 & 0 & 0 & 0 & 0 & 0 & 0 & 0 \\
   0\vphantom{(0')} & 1 & 1 & 1 & 0 & 0 & 0 & 0 & 0 & 0 & 0 & 0 \\
   1\vphantom{(0')} & 0 & 0 & 0 & 0 & 0 & 0 & 0 & 0 & 0 & 0 & 0 \\
   0\vphantom{(0')} & 1 & 1 & 1 & 0 & 0 & 0 & 0 & 0 & 0 & 0 & 0 \\
   1\vphantom{(0')} & 0 & 0 & 0 & 0 & 0 & 0 & 0 & 0 & 0 & 0 & 0 \\
   1\vphantom{(0')} & 0 & 0 & 0 & 0 & 0 & 0 & 0 & 0 & 0 & 0 & 0 \\
   1\vphantom{(0')} & 0 & 0 & 0 & 1 & 0 & 0 & 0 & 0 & 0 & 0 & 0 \\
  \end{smallmatrix}\right)}_{\mS},
  \;\;
  \underbrace{\left(\begin{smallmatrix}
   0.6\vphantom{(0')} \\
   0.9\vphantom{(0')} \\
   0.3\vphantom{(0')} \\
   0.1\vphantom{(0')} \\
   0.9\vphantom{(0')} \\
   0.6\vphantom{(0')} \\
   0.3\vphantom{(0')} \\
   0.5\vphantom{(0')} \\
   0.1\vphantom{(0')} \\
   0.4\vphantom{(0')} \\
   0.5\vphantom{(0')} \\
   1.0\vphantom{(0')} \\
  \end{smallmatrix}\right)}_{\vw},
  \;\;
  \underbrace{\left(\begin{smallmatrix}
   1\vphantom{(0')}  & 1  & 1  & 1  & 0  & 0  & 0  & 0  & 0  & 0  & 0  & 0 \\
   0\vphantom{(0')}  & 0  & 0  & 0  & 1  & 1  & 0  & 0  & 0  & 0  & 0  & 0 \\
   0\vphantom{(0')}  & 0  & 0  & 0  & 0  & 0  & 1  & 1  & 0  & 0  & 0  & 0 \\
   0\vphantom{(0')}  & 0  & 0  & 0  & 0  & 0  & 0  & 0  & 1  & 1  & 0  & 0 \\
   0\vphantom{(0')}  & 0  & 0  & 0  & 0  & 0  & 0  & 0  & 0  & 0  & 1  & 0 \\
   0\vphantom{(0')}  & 0  & 0  & 0  & 0  & 0  & 0  & 0  & 0  & 0  & 0  & 1 \\
   1\vphantom{(0')}  & 0  & 0  & 0  & 1  & 0  & 1  & 0  & 1  & 0  & 0  & 0 \\
   0\vphantom{(0')}  & 1  & 0  & 0  & 0  & 1  & 0  & 0  & 0  & 0  & 0  & 1 \\
   0\vphantom{(0')}  & 0  & 1  & 0  & 0  & 0  & 0  & 1  & 0  & 0  & 0  & 0 \\
   0\vphantom{(0')}  & 0  & 0  & 1  & 0  & 0  & 0  & 0  & 0  & 1  & 0  & 0 \\
   0\vphantom{(0')}  & 0  & 0  & 0  & 0  & 0  & 0  & 0  & 0  & 0  & 1  & 0 \\
 \end{smallmatrix}\right)}_{\mC}
 $

 \caption{A small sample problem and the data for the QP formulation.}
 \label{fig:netalign-qp-data}
 \end{figure*}

Using these definitions, the network alignment problem is an integer quadratic program (QP)
\begin{equation} \tag{NAQP} \label{eq:naqp}
 \MAXone{\vx}{\alpha \vw^T \vx + \beta/2 \vx^T \mS \vx}{%
 \mC \vx \le \vI_{n+m}, \qquad x_{ii'} \in \{0, 1\}}
\end{equation}
where $\alpha$ and $\beta$ are arbitrarily chosen nonnegative constants
that define the trade-off between the similarity and overlap objectives.
When $\alpha=0$ and $\beta=1$ then the
program solves a special case we
call the overlap graph matching problem or
the pure overlap problem.
When $\alpha=1$ and $\beta=0$
it solves the maximum weight matching problem.
Figure~\ref{fig:netalign-qp-data} shows an example of the network
alignment problem and an explicit construction of the matrices $\mS$
and $\mC$.

\section{Applications}\label{sec:relation2other-works}
Network alignment is deeply intertwined with many classical
computational problems such as graph isomorphism, quadratic assignment,
maximum common subgraph, and maximum clique. For a survey of these
connections see \citet{conte2004-graph-matching}.
In this section we briefly highlight the key applications
of network alignment that appear in
pattern recognition, ontology alignment and bioinformatics.

\subsection{Pattern recognition}
Network alignment for pattern recognition involves
identifying a small model graph within a
large scene graph.  The model graph typically represents
the desired pattern -- a rooftop, a face, a person --
and the scene graph describes the entire space -- possibly
a picture.  The assumptions are often that the data are noisy
and the goal is not an exact subgraph isomorphism.
Again, \citet{conte2004-graph-matching} is a good
starting point to explore this literature.
Recent work includes trying to \emph{learn}
good scores $\vw$ to avoid using the quadratic
matching formulations~\cite{Caetano-2009-graph-matching}.

\subsection{Ontology matching}

An ontology is a set of statements, which connect subjects to
objects with verbs.  An elementary example is an ontology
describing the authors of this paper:
\begin{center}
\begin{tabular}{rcl}
subject & verb & object \\
David Gleich & \emph{wrote} & Sparse Network Alignment \\
Mohsen Bayati & \emph{wrote} & Sparse Network Alignment \\
Sparse Network Alignment & \emph{is an} & Academic Manuscript \\
Academic Manuscript & \emph{is a} & Paper\\
\end{tabular}
\end{center}
and so on.  An ontology is a flexible data description
format, and a fundamental problem is how to align two ontologies
about the same data.  Suppose that Citeseer and DBLP expose their
networks of papers as an ontology, the
problem of ontology alignment is to figure out the correspondence
between Citeseer papers and DBLP papers. This problem has been studied extensively.
See \citet{hu2008-falcon-ao}, \citet{ehrig2004-qom}, \citet{hu2005-gmo}, and
\citet{blondel2004-graph-similarity} for a few different approaches
to these problems.  All of these approaches utilize some heuristic
approach for a network alignment problem.

\subsection{Finding common pathways in biological networks}

Network alignment is becoming a small industry
within bioinformatics.  Broadly speaking, this emergence
is due to the rapid increase in high quality data about
protein interactions.  A protein-protein interaction (\acro{PPI})
graph has proteins as vertices and edges that connect proteins
known to interact.  Suppose that $\graph{A}$ and $\graph{B}$
are two \acro{PPI} networks, and we compute an alignment between
them.  The alignment produces a one-to-one mapping between proteins
in $\graph{A}$ and proteins in $\graph{B}$.  If the proteins
are from two different species, then the alignment hints
at similar functions for the two proteins, or two groups of proteins.
Alternatively, we may know information about proteins
in $\graph{A}$.  An alignment with $\graph{B}$ suggests
what information about $\graph{A}$ might apply to the proteins
in $\graph{B}$.

Due to the wide interest in this problem, several tools have been developed for aligning protein-protein interaction networks.
These include NetAlign~\cite{liang2006-netalign},
Gr\ae mlin~\cite{flannick2006-graemlin,%
flannick2008-graemlin2} and IsoRank~\cite{singh2007-matching-topology,%
singh2008-isorank-multi}.
Some of these tools have extensions for aligning more than two networks,
but we focus on the two network case here.  We
review the IsoRank algorithm in detail in Section~\ref{sec:isorank}.
Alternative approaches are proposed by \citet{berg2006-bayesian-alignment} and
\citet{kuchaiev2009-topological-alignment}.

\section{Existing algorithms for network alignment}\label{sec:existing-algorithms}

In this section we review existing algorithms that produce good solutions for
the network alignment. 

\subsection{IsoRank Algorithm}\label{sec:isorank}

\citet{singh2007-matching-topology} proposed IsoRank
to  approximately solve \ref{eq:naqp} when
$\graph{L}$ is a complete bipartite graph.
In this section, we present IsoRank and our variation SpaIsoRank,
which is more efficient when $\graph{L}$ is sparse.

The main idea of IsoRank algorithm is to approximate the objective
of \ref{eq:naqp} without direct concern for the matching constraints.
Let $\mA$ and $\mB$
be the adjacency matrices for graphs $A$ and $B$, and also let $\mD_A$
and $\mD_B$ be the diagonal matrices of their degrees, respectively.  IsoRank solves
for the matrix $\mZ$
that satisfies:
\[ \gamma \underbrace{\mA^T \mD_A}_{\mP^T} \mZ \underbrace{\mD_B \mB}_{\mQ} + (1-\gamma) \mW = \mZ. \]
Here $W_{i,j} = w_{i,j}$ is the weight function on edges in $L$ represented as a matrix
when $L$ is the complete bipartite graph.
The intuition is that each entry $Z_{i,i'}$ is a real number based on a weighted average of all neighboring
values $Z_{j,j'}$ where $(i,j)\in E_A$ and $(i',j')\in E_B$.
With this heuristic solution $\mZ$, they compute a binary solution $\mX$ by solving
a maximum weight matching problem where the weights are from $\mZ$.
We discuss rounding schemes in more detail in Section~\ref{sec:rounding}.

Now we discuss our extensions of IsoRank algorithm for the case when $L$ is sparse.  This extension rests on the
Kronecker product.
Recall that the Kronecker product of an $m \times n$ matrix
$\mA$ and another matrix $\mB$ is defined by
\begin{equation}
 \mA \otimes \mB =
   \pmat{ A_{11} \mB & \cdots & A_{1n} \mB\\
          \vdots & \ddots & \vdots \\
          A_{m1} \mB & \cdots & A_{mn} \mB}
\end{equation}
When $L$ is the complete bipartite graph, then the matrix
indicating potential overlaps, or squares, is
$\mS = \mB \kron \mA$ (up to an arbitrary permutation
based on the edge order $\order_L$).  Moreover,
the mixed product property states that
\[ \mP^T \mZ \mQ = (\mQ \kron \mP) \tvec(\mZ), \]
where $\tvec(\mZ)$ is a column-wise vector representation of
a matrix:
\[ \tvec(\mZ) = \bmat{\mZ \ve_1 \\ \mZ \ve_2 \\ \vdots \\ \mZ \ve_n}. \]
Note that $\mQ \kron \mP = \diag[\mS \vI_{|E_L|}] (\mB \kron \mA)$.

Thus, the following PageRank problem is equivalent to IsoRank
when $L$ is the complete bipartite graph, but handles
sparse $L$ as well:
\[ \gamma \mD_S \mS^T \vz + (1-\gamma) \vw = \vz. \]
We compute $\vz$ using a standard algorithm for PageRank.
At each iteration, we employ one of two rounding schemes to
produce a matching.  The first just uses the vector $\vz\itn{k}$
as the weight on each edge in $L$ and solves a bipartite max-weight
matching problem.  The second uses the vector $\alpha \vw + (\beta/2) \mS \vz\itn{k}$
as the weights on $\mL$ and more closely mirrors the original
objective function.
\begin{inlinealgorithm}{SpaIsoRank}{$\mS$, $L$,
     damping (nonnegative) parameter $\gamma < 1$, $\eps$, n$_{\textrm{iter}}$,
     rounding\_type $\in \{ 1, 2 \}$ }
\begin{lstlisting}[language=Matlab,
  numbers=left,numberstyle=\scriptsize,xleftmargin=5mm,
  mathescape=true,columns=flexible]
$\vv = \vw/(\vI_{|E_L|}^T \vw)$
$\vd = \mS \vI_{|E_L|}$, $\mP = \diag[\vd]^{-1} \mS$
$\vz\itn{0} = \vv$, $\delta=\eps+1$
for $k=1$ to n$_{\textrm{iter}}$ unless $\delta < \eps$
  $\vz\itn{k} = \gamma \mP^T \vz\itn{k-1} + (1-\gamma) \vv$
  $\delta = \normof{\vz\itn{k}-\vz\itn{k-1}}$
  if rounding_type is 1
    $\vx\itn{k} = $ bipartite_match(L, $\vz\itn{k}$)
    $\text{obj}\itn{k} = $ objective($\vx\itn{k}$)
  else if rounding_type is 2
    $\vx\itn{k} = $ bipartite_match(L, $\alpha \vw + (\beta/2) \mS \vz\itn{k}$)
    $\text{obj}\itn{k} = $ objective($\vx\itn{k}$)
  end
end
return $\vx\itn{k}$ with the highest value of $\text{obj}\itn{k}$
\end{lstlisting}
\end{inlinealgorithm}

\SubSection{Linear Program Formulations}

We now review a series of linear programming (LP) relaxations for
network alignment. These ideas originated in mixed integer
translations of the quadratic assignment problem \cite{lawler1963-qap},
and subsequent tightened versions by that were originally described
by \citet{Frieze1983-qap} and \citet{Adams1994-improved-qap}.
The adaptation to network alignment appeared in
\citet{klau2009-network-alignment}.

In the first relaxation,\citet{lawler1963-qap},
converted \ref{eq:naqp} into a mixed integer
linear program.  To do so,  replace each
product $x_{ii'} x_{jj'}$ with a new variable $y_{ii',jj'}$,
and add constraints $y_{ii',jj'} \le x_{ii'}$, and $y_{ii',jj'} \le
x_{jj'}$. These constraints enforce $y_{ii',jj'} \le x_{ii'}
x_{jj'}$ when $x_{ii'}$ and $x_{jj'}$ are binary.
We also add symmetry constraints $y_{ii',jj'}=y_{jj',ii'}$. Notice that with the
symmetry constraints the constraints $y_{ii',jj'} \le x_{jj'}$ can be dropped.

Before writing the new integer program, let us define $\mat{Y}_{\mS}$ to be a matrix with the same dimension as $\mS$ where
\[
\mat{Y}_{\mS}[ii',jj']=
\left\{
\begin{array}{ll}
y_{ii',jj'}&\textrm{ if }\mS[ii',jj']=1\\
0&\textrm{ Otherwise.}
\end{array}
\right.
\]
Thus, we arrive at:
\begin{equation} \tag{NAILP} \label{eq:nailp}
 \MAXthree{\vx,\vy}{\alpha \vw^T \vx + \frac{\beta}{2}\sum_{ii'}\sum_{ii'\Box jj'}y_{ii',jj'}}{%
    \mC \vx \le \vI_{n+m}, \qquad x_{ii'} \in \{0, 1\},}{%
    y_{ii',jj'} \le x_{ii'} \qquad \text{ for all } ii'\Box jj',}{%
    \mat{Y}_{\mS}=\mat{Y}_{\mS}^T}
\end{equation}
as a mixed-integer linear program to solve the network alignment problem

In contrast with the quadratic program, we can relax the binary
constraint on \ref{eq:nailp} and get an efficient algorithm.
After we write $\sum_{ii'}\sum_{ii'\Box jj'}y_{ii',jj'}$ as
$\mS\bullet \mat{Y}_{\mS}$, the relaxed program is
\begin{equation} \tag{NARLP} \label{eq:narlp}
 \MAXthree{\vx,\vy}{\alpha \vw^T \vx + \frac{\beta}{2}\mS\bullet \mat{Y}_{\mS}}{%
    \mC \vx \le \vI_{n+m}, \qquad x_{ii'} \in [0, 1],}{%
    y_{ii',jj'} \le x_{ii'} \qquad \text{ for all } ii'\Box jj',}{%
    \mat{Y}_{\mS}=\mat{Y}_{\mS}^T}.
    \end{equation}
It admits a polynomial-time solution with an appropriate linear program
solver.

\begin{remark} The relaxation \ref{eq:narlp} is
advantageous because it yields
an upper bound on the objective value of the network alignment problem.
Furthermore, solving \ref{eq:narlp} with $\alpha=0,\beta=1$
allows us to get an upper bound on the maximum possible
overlap between two networks.
\end{remark}

\SubSection{Klau's Iterative Matching Relaxation}\label{subsec:Klau}

\citet{klau2009-network-alignment} constructed an iterative
algorithm to approximate \ref{eq:naqp}.  The key components
of this algorithm are a tighter LP relaxation of
\ref{eq:naqp} and the Lagrangian decomposition of the
symmetry constraints.  We first explain the Lagrangian
decomposition for \ref{eq:narlp} and then show the
tightened LP. In the Lagrangian decomposition,
we drop all the symmetry constraints $\mat{Y}_{\mS}=\mat{Y}_{\mS}^T$
by adding penalty terms of
the form $u_{ii',jj'}(y_{ii',jj'}-y_{jj',ii'})$.  Here
$u_{ii',jj'}$'s are Lagrange multipliers, a set of $n^2 - n$
new variables. Following
this idea, we arrive at
\begin{equation} \tag{NALLP} \label{eq:nallp}
 \MAXtwo{\vx,\vy}{\alpha \vw^T \vx + \frac{\beta}{2}\mS\bullet \mat{Y}_{\mS}+\mat{U}_{\mS}\bullet(\mat{Y}_{\mS}-\mat{Y}_{\mS}^T)}{%
    \mC \vx \le \vI_{n+m}, \qquad x_{ii'} \in [0, 1],}{%
    y_{ii',jj'} \le x_{ii'} \qquad \text{ for all } ii'\Box jj'.}
\end{equation}
When $\mat{Y}_{\mS}=\mat{Y}_{\mS}^T$, the two linear
programs \ref{eq:narlp} and \ref{eq:nallp} are equivalent.
Therefore, for any fixed $\mat{U}_{\mS}$ the optimum solution of
\ref{eq:nallp} is an upper bound for the objective of
\ref{eq:narlp}, which is itself an upper bound for the network
alignment problem. Standard Lagrangian theory dictates that
with the optimal Lagrange multipliers $\mat{U}_{\mS}$,
the two LP's have the same optimum.  The advantage of
using \ref{eq:nallp} is that the solution is integral
for any fixed $\mU_{\mS}$, and moreover, we can compute
it by solving a max-weight matching problem.
Let us explain why that happens.
For a fixed $\mU$, note that the objective decouples
between $\vx$ and $\vy$:
\[ \begin{aligned} & \alpha \sum_{ii'} w_{ii'} x_{ii'} + \frac{\beta}{2} \sum_{ii' \Box jj'} y_{ii',jj'} + \sum_{ii \Box jj'} u_{ii',jj'}(y_{ii',jj'} - y_{jj',ii'}) \\
& \quad  = \alpha \sum_{ii'} w_{ii'} x_{ii'} +  \sum_{ii' \Box jj'} y_{ii',jj'} (\frac{\beta}{2} + u_{ii',jj'} - u_{jj',ii'}).
\end{aligned} \]
Because $y_{ii',jj'} \le x_{ii'}$, the optimum is
\[ y_{ii',jj'} = \begin{cases} 0 &  \frac{\beta}{2}+u_{ii',jj'} - u_{jj',ii'} < 0 \\ x_{ii'} & \text{otherwise}. \end{cases} \]
Therefore, let
\[ \bar{w}_{ii'} = \alpha w_{ii'} +
    \sum_{ii' \Box jj'}\max \{0, \frac{\beta}{2}+u_{ii',jj'} - u_{jj',ii'}\}.
\]
Then the solution of \ref{eq:nallp} can be found by solving the following max-weight-matching problem:
\[
\MAXone{\vx}{\bar{\vw}^T \vx}{\mC \vx \le \vI_{n+m}, \qquad x_{ii'} \in \{ 0, 1 \}}.
\]
Thus, for any fixed Lagrange multipliers, we can solve \ref{eq:nallp} as a single max-weight
matching problem.  In effect, we have grouped the objective
function of \ref{eq:nallp} into pieces where $y_{ii',jj'}$
is completely determined by $x_{ii'}$.  Futhernote that
we if $ u_{ii',jj'} - u_{jj',ii'} = 0$, then we get an especially
simple means of upper-bounding the overlap with a single
max-weight matching.
%
%
%

While these relaxations give upper bounds on the objective, there is
often a large gap between the upper bound and the integer
solution.  \citet{Frieze1983-qap} and \citet{Adams1994-improved-qap}
propose tightened LPs for the quadratic assignment problem.
Klau's algorithm adapts these improvements to the network
alignment problem.  Notice
that in both \ref{eq:nailp} and \ref{eq:narlp}
\begin{equation}
\label{eq:yii'jj'-is-matching}
\sum_j y_{ii',jj'}\leq \sum_j x_{jj'}\leq 1,
\qquad
\sum_{j'} y_{ii',jj'}\leq \sum_{j'} x_{jj'}\leq 1.
\end{equation}
for any fixed $ii'$.
This means that row $ii'$ of $\mat{Y}_{\mS}$ (denoted by
$\mat{Y}_{\mS}[ii',:]$) should satisfy the matching constraint
$\mC\big( \mat{Y}_{\mS}[ii',:]\big)^T\le \vI_{n+m}$. However, when the
symmetry constraints are removed, the inequalities \eqref{eq:yii'jj'-is-matching} may
be violated.  The tightened LP re-adds these constraints:
\begin{equation} \tag{NATLP} \label{eq:natlp}
 \MAXthree{\vx,\vy}%
    {\alpha \vw^T \vx + \frac{\beta}{2}\mS\bullet\mat{Y}_{\mS} + %
       \mat{U}_{\mS}\bullet(\mat{Y}_{\mS}-\mat{Y}_{\mS}^T)}%
    {\mC \vx \le \vI_{n+m}, \qquad x_{ii'} \in [0, 1],}%
    {y_{ii',jj'} \le x_{ii'} \qquad \text{ for all } ii'\Box jj',}%
    {\mC\big( \mat{Y}_{\mS}[ii',:]\big)^T\le \vI_{n+m}
       \qquad \text{ for all } ii' }
\end{equation}
This tighter LP still can be solved using a MWM algorithm. Like in \ref{eq:nallp}, $y_{ii',jj'}$ can be grouped by $ii'$. Now the term
\[
\max_{\vy}\sum_{ii'\Box jj'}y_{ii',jj'}(\frac{\beta}{2}+u_{ii',jj'}-u_{jj',ii'})
\]
equals $x_{ii'}$ times the solution of a \emph{small} MWM problem
\begin{equation} \label{eq:smallMWM}
 \MAXone{\mat{Y}_{\mS[ii',:]}}%
   {\mat{Y}_{\mS}[ii',:]
        \left[\frac{\beta}{2}\vI_{|E_L|}^T +
              \mat{U}_{\mS}[ii',:]-\mat{U}_{\mS}^T[ii',:]
        \right]^T}%
   {\mC\big( \mat{Y}_{\mS}[ii',:]\big)^T\le \vI_{n+m}, \qquad y_{ii',jj'} \in \{ 0, 1 \}}.
\end{equation}


Klau's final algorithm is an iterative procedure that uses
a sub-gradient algorithm to optimize $\mU_{\mS}$ (step 8 below).
Each step of the sub-gradient method involves solving \ref{eq:natlp}
for a new $\mU_{\mS}$.  These sub-problems are solved by solving
$|E_L|$ small max-weight matching problems for \eqref{eq:smallMWM}
and then a single large max-weight matching problem for \ref{eq:natlp}.
Because each iteration of this algorithm calls many max-weight
matching functions, we call this algorithm the \emph{matching relaxation} or
MR for short.

To state the algorithm compactly, we need a small
bit of new notation.
First, let $a \le b$.  Define
  \[ \bound_{a,b} z \equiv \min(b,\max(a,z))
     = \begin{cases} a & z < a \\ z & a \le z \le b \\ b & z > b \end{cases} \]
  and let both $\bound_{a,b} \vx$ and $\bound_{a,b} \mA$ be defined element-wise.
Also define
\[ \vd,\mat{S_L} = \text{maxrowmatch}(\mS,L) \]
where each entry $ii'$ in $\vd$ is the result of a MWM on all the other
edges in row $ii'$ of $\mS$, with weights from the corresponding
entries of $\mS$.  Written formally, $d_{ii'} = $
bipartite\_match($\{ jj' : S[ii',jj']=1) \})$. The matrix
$\mat{S_L}$ has a 1 for any edge used in the optimal solution of the
bipartite matching problem for a row.
\begin{inlinealgorithm}{NetAlignMR}{$\mS$, $\vw$,
             nonnegative damping parameters $\gamma\leq 1$, n$_{\textrm{iter}}$, mstep, $\alpha$, $\beta$}
\begin{lstlisting}[language=Matlab,
  numbers=left,numberstyle=\scriptsize,xleftmargin=5mm,
  mathescape=true,columns=flexible]
$\mU\itn{0} = 0$
for $k=1$ to n$_{\textrm{iter}}$
  $\vd,\mat{S_L}$ = maxrowmatch($(\beta/2) \mS + \mU - \mU^T, L$)
  $\bar{\vw}\itn{k} = \alpha \vw + \vd$
  $\vx\itn{k}$ = bipartite_match($L$,$\bar{\vw}\itn{k}$)
  $\text{obj}\itn{k} = \alpha {\vx\itn{k}}^T \vw + \beta/2 {\vx\itn{k}}^T \mS \vx\itn{k}$ % \objective_function($\vx\itn{k}$)
  $\text{upper}\itn{k} = \bar{\vw}^{(k)^T} \vx\itn{k}$
  $\mF = \mU\itn{k-1} - \gamma \mX\itn{k} \text{\bfseries triu}(\mat{S_L}) + \gamma \text{\bfseries tril}(\mat{S_L})^T \mX\itn{k}$
  $\mU\itn{k} = {\displaystyle \bound_{-0.5,0.5} } \mF $
  if $\text{upper}\itn{k}$ has not changed in mstep iterations,
    set $\gamma = \gamma/2$
  end
end
return $\vx\itn{k}$ with the largest value of $\text{obj}\itn{k}$
\end{lstlisting}
\end{inlinealgorithm}

Here $\text{\bfseries triu}(\mat{S_L})$ ($\text{\bfseries tril}(\mat{S_L})$)
represents the upper (lower) triangular part of $\mat{S_L}$.
We reduce the sub-gradient step-length $\gamma$ on a schedule
determined by the change in the upper-bound.

%
%

\Section{Our result: Two message passing algorithms} \label{sec:bp}

In this section, we introduce two message passing algorithms for
network alignment. Message passing has been remarkably successful
in coding theory \cite{Gal63}, artificial intelligence \cite{Pea88}, solving constraint satisfaction problems \cite{MeZ02}, structural biology \cite{YaW02}, computer vision \cite{TaF03}, data clustering \cite{frey2007-affinity}, and compressed sensing \cite{DMM09}.

More precisely, we will use a Belief Propagation (BP) approach. In
general, BP works by iteratively making local and greedy decisions.
Decisions are updated by passing messages between neighboring entities (nodes of the graph).

In what follows, we provide a quick overview of related BP approaches for the matching problems, next we derive a BP based algorithm that passes messages
along the edges of graph $\graph{L}$ and also among the squares
(Section~\ref{sec:netalignmp}).  These messages have an intuitive representation, which we present in Section~\ref{subsec:bp-intuition}.
Next, we state a matrix-based version of the same algorithm
(Section~\ref{subsec:bp-matrix-form}) in order to elucidate the data
organization and computation.  Finally, we conclude by developing a
more technical message passing algorithm that includes additional
constraints from \ref{eq:natlp} (Section~\ref{sec:netalignmp++}).

\SubSection{Related work on BP and graph matching.} BP approaches have been shown to correctly find the optimum solution for a variety of optimization problems including maximum-weight matching \cite{BSS05}, \cite{SMW11}, and \cite{BBCZ07} and our algorithm for the network alignment problem is inspired by that approach. However, the matching problem studied in \cite{BSS05}, \cite{SMW11}, and \cite{BBCZ07} is a very special case of the network alignment problem (when $\beta=0$ in \ref{eq:naqp}) that can be solved exactly in polynomial time. The quadratic term that appears when $\beta\neq 0$ is NP hard to maximize and requires special treatment which we carry by defining a new factor graph on squares and edges of the graph $\graph{L}$. Recently, and independently from our work, \cite{BBMTWZ10} introduced a completely different BP approach for aligning graphs in biology using a maximum clique representation of the problem.


%
%
\SubSection{A factor graph representation}\label{subsec:factor-graph}

To use BP, it is standard to define a
probability distribution on the space of
all matchings in $L$ that assigns the highest probability to the
matching that maximizes \ref{eq:naqp}.
This matching is also called the \emph{maximum a posteriori assignment} (MAP).
We begin with this construction.

Let $V_A = \{1,\ldots,n\}$ and $V_B = \{1',\ldots,m\}$. For any square
formed by the two edges $ii'$ and $jj'$ of $E_L$, we create a new vertex $ii'jj'$, and denote
the set of all such vertices by $V_S$, i.e.
\[
V_S = \left\{ii'jj'|~ ii' \textrm{ and } jj' \textrm{ form a square }\right\}\,.
\]
Now, we assign a binary variable $x_{ii'}$ to each edge $ii'\in E_L$ and a binary variable $x_{ii'jj'}$ for each square $ii'jj' \in V_S$.
We also use notation $\vx_D$ for any subset $D\subset E_L\cup V_S $ to denote the
vector $[x_d]_{d\in D}$.  The set of neighbors of a node $v$ in a graph $G$ is denoted by $\pa v$.

\begin{figure}
\centering
      \begin{tabular}{c@{\hspace*{4em}}c}
       \includegraphics{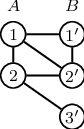} & \includegraphics{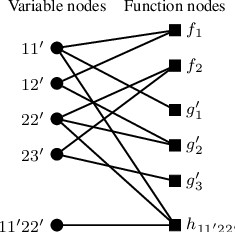} \\
      (a) & (b) \\
      \end{tabular}
      \caption{\normalsize The graph (b) is the factor-graph representation of the network alignment problem in (a).}
      \label{fig:fact-graph}
\end{figure}
%

%

Next, we define a new graph (factor graph) that has the following two types of nodes:
\begin{itemize}
\item[(i)] \emph{Variable nodes:} $|E_L|+|V_S|$ nodes, one for each element of $E_L$ and $V_S$.
The binary variables assigned to these nodes are denoted by $(\vx_{E_L},\vx_S) \in \{0,1\}^{|E_L|+|S|}$.

\item[(ii)] \emph{Function nodes:} $|V_A|+|V_B|+|S|$ nodes of two types. One type is
for enforcing the integer constraints. That is for each vertex $i\in V_A$ ($i'\in V_B$)
we define a function node $f_{i}:\{0,1\}^{|E_L|+|S|}\to \mb{R}$ ($g_{i'}:\{0,1\}^{|E_L|+|S|}\to \mb{R}$) by:
\begin{align*}
f_{i}\big(\vx_{\pa f_i}\big)
  & = \begin{cases} 1 & \sum_{ii'\in E_L} x_{ii'}\leq 1 \\ 0 & \text{otherwise} \end{cases}
   && \text{ for all } i
   \\
g_{i'}\big(\vx_{\pa g_{i'}}\big)
  & = \begin{cases} 1 & \sum_{ii'\in E_L} x_{ii'}\leq 1 \\ 0 & \text{otherwise} \end{cases} && \text{ for all } i'\,.
\end{align*}
The neighbor operation used to define the left-hand vector $\vx_{\pa f_i}$
is implicitly defined by the set of variables used on the right-hand
side of the equation.
In words, the function node $f_{i}$ ($g_{i'}$) enforces the matching
constraint at $i$ ($i'$)

Another type of function nodes check the validity of squares. For each square $ii' \Box jj'$ define a function node
$h_{ii'jj'}:\{0,1\}^{|E_L|+|S|}\to \mb{R}$:
\begin{align*}
h_{ii'jj'}\big(\vx_{\pa h_{ii'jj'}}\big)
& = \begin{cases} 1 & x_{ii'jj'}=x_{ii'}x_{jj'} \\ 0 & \text{otherwise} \end{cases} && \text{ for all } (ii',jj') \in V_S\,.
\end{align*}
In other words, $h_{ii'jj'}$ guarantees that $x_{ii'jj'}=1$ if and only if $x_{ii'}= x_{jj'}=1$.
\end{itemize}

The edges of the factor graph are simply connecting each function node to the variable nodes it acts on. For example each $f_i$ is connected to all variable nodes $ii'\in E_L$ and each $h_{ii'jj'}$ is connected to $ii'$, $jj'$ and $ii'jj'$ in $E_L\cup V_S$. Therefore the factor graph is bipartite.

Figure \ref{fig:fact-graph} shows an example of a graph pair $\graph{A}, \graph{B}$ and their factor-graph representation as described above.

Now define the following probability distribution
\begin{align}
\label{eq:gibbs-prob}
p(\vx_{L},\vx_S)=\frac{1}{Z}\left[
\dsp\prod_{i=1}^nf_{i}(\vx_{\pa f_i})\dsp\prod_{j=1}^{m}g_{j}(\vx_{\pa g_j})\dsp\prod_{ijrs\in V_S}h_{ijrs}(\vx_{\pa h_{ijrs}})\right]
e^{\alpha\vw^T \vx_{L}+\frac{\beta}{2} \vI_{|S|}^T\vx_S}
\end{align}
where $Z$ is just a normalization term to make $p(\vx_{L},\vx_S)$ a probability distribution.
In particular,
\[
Z\equiv\sum_{(\vx_{L},\vx_S)\in \{0,1\}^{|E_L|+|S|}}\left[
\dsp\prod_{i=1}^nf_{i}(\vx_{\pa f_i})\dsp\prod_{j=1}^{m}g_{j}(\vx_{\pa g_j})\dsp\prod_{ijrs\in V_S}h_{ijrs}(\vx_{\pa h_{ijrs}})\right]
e^{\alpha\vw^T \vx_{L}+\frac{\beta}{2} \vI_{|S|}^T\vx_S}\,.
\]
Note that, there is a 1-1 correspondence between the feasible solutions of
\ref{eq:naqp} and support of the probability distribution \eqref{eq:gibbs-prob}.
The following lemma formalizes this observation.
\begin{lemma}\label{lem:1-1corresp}
For any $(\vx_{L},\vx_S) \in \{0,1\}^{|E_L|+|V_S|}$ with non-zero probability, the vector $\vx_{L}$ satisfies the constraints of the integer program \ref{eq:naqp}. Conversely, any feasible solution $\vx_{L}$ to \ref{eq:naqp} has a unique counterpart $(\vx_{L},\vx_S)$ with non-zero probability $p(\vx_{L},\vx_S)=e^{\alpha \vw^T \vx + (\beta/2) \vI_{|S|}^T\vx_S}$.
\end{lemma}
\begin{proof}
Any $(\vx_{L},\vx_S) \in \{0,1\}^{|E_L|+|V_S|}$ with non-zero probability should satisfy the conditions dictated by function nodes $f, g, h$ which translates to $\vx_{L},\vx_S$ being a feasible solution to \ref{eq:naqp}. Conversely, for any feasible solution to \ref{eq:naqp} the values of function nodes $f,g,h$ are equal to $1$ and hence the probability is non-zero.
\end{proof}

Moreover, any pair with maximum probability is an optimum solution to \ref{eq:naqp}.
\begin{lemma}\label{cor:MAP-2-OPT}
The vector $(\vx_{L}^*,\vx_S^*)$ is equal to $\arg\max_{\vx_{L},\vx_S}{p(\vx_{L},\vx_S)}$
if and only if $\vx_{L}^*$ is the optimum solution to \ref{eq:naqp} and $\vx_S^*$ is the vector of squares generated by it.
\end{lemma}
\begin{proof}
Proof immediately follows from Lemma \ref{lem:1-1corresp}.
\end{proof}
Using Lemma \ref{cor:MAP-2-OPT}, it is known that a variant of BP algorithm (\emph{max-product} or \emph{min-sum}) can be used to find an approximate solution to $\ref{eq:naqp}$ \cite{Mezard-Montanari-book}. In this paper we use the notion BP to refer to this variant.
%
%
\subsection{The message passing algorithm}\label{sec:netalignmp}
The standard BP messages for finding the optimum solution $\arg\max_{\vx_{L},\vx_S}{p(\vx_{L},\vx_S)}$ are vectors of numbers. However, for our problem we show that the information contained in these vector messages can be compressed to a real number. Therefore,
we can obtain a simple algorithm with a smaller running time that will be presented next.  For completeness, we provide the derivation of this simplified version from the standard BP in Appendix \ref{sec:BP-eqs}. However, in Section \ref{subsec:bp-intuition} we provide a more intuitive description of the algorithm.

\begin{inlinealgorithm}{NetAlignMP}{$\alpha, \beta$, the set of squares $V_S$, and
  the weighted bipartite graph $\graph{L} = (V_A \cup V_B,E_L)$,
  and a damping parameter $\gamma$.}

\begin{enumerate}
\item At times $t=0,1,\ldots$, each edge $ii'$ sends two messages of the form $m_{ii'\to f_i}^{(t)}$ and $m_{ii'\to g_{i'}}^{(t)}$ and also sends one message of the form $m_{ii'\to h_{ii'jj'}}^{(t)}$ for any square $ii'\Box jj'$.

\item Initialize messages to 0.

\item For $t \geq 1$, the messages in iteration $t$ are obtained from the
messages in iteration $t-1$. In particular for all $ii'\in E_L$
\begin{equation}
 \label{eq:bp1}
 \begin{aligned}
m_{ii'\to f_i}^{(t)} = & \alpha w_{ii'} -
  \left(\max_{k\neq i}\left[m_{ki'\to g_{i'}}^{(t-1)}\right]\right)_+ \;\;  + \;\;
  \sum_{\mathclap{jj' : ii'\Box jj'}}
                \bigl[ (\frac{\beta}{2}+m_{jj'\to h_{ii'jj'}}^{(t-1)})_+ -
                       (m_{jj'\to h_{ii'jj'}}^{(t-1)})_+
                \bigr].
\end{aligned}
\end{equation}
Here, notation $(x)_+$ represents $\max(0,x)$. The update rule for $m_{ii'\to g_{i'}}^{(t)}$ is similar to the update rule for $m_{ii'\to f_i}^{(t)}$ and
\begin{equation}\label{eq:bp2}
\begin{aligned}
m_{ii'\to h_{ii'jj'}}^{(t)} =
  & \alpha w_{ii'} +
       \sum_{ii'kk'\neq ii'jj'}
             \bigl[ (\frac{\beta}{2} + m_{kk'\to h_{ii'kk'}}^{(t-1)})_+ -
                    (m_{kk'\to h_{ii'kk'}}^{(t-1)})_+\bigr] \\
  & \qquad - \left(\max_{k\neq i}\left[m_{ki'\to g_{i'}}^{(t-1)}\right]\right)_+
           - \left(\max_{k'\neq i'}\left[m_{ik'\to f_i}^{(t-1)}\right]\right)_+.
\end{aligned}
\end{equation}



\item Apply damping on the message updates.  (See possibilities in
    Section~\ref{sec:BP-Convergence}.)
\item Round the solution (see possibilites in Section~\ref{sec:rounding})
  and compute the objective function on the rounded messages.

\item Repeat (3)-(5) for a fixed number of iterations
unless the messages stop changing.
\end{enumerate}
\normalsize {\small \sffamily OUTPUT} the rounded solution with the best objective value.
\end{inlinealgorithm}

%
%
\subsection{The intuition behind \namp}\label{subsec:bp-intuition}


\namp exploits the fact that the constraints of \ref{eq:naqp} are
local.  Suppose each edge of the graph $\graph{L}$ is an agent and
each agent can talk to its neighbors. First observe that together,
the agents can verify the feasibility of any solution to
\ref{eq:naqp}. The next step is to note that they can also calculate
the merit of each solution ($\alpha \vw^T\vx+\beta/2 \vx^T\mS\vx$)
locally.

Based on the above intuition, each agent should communicate to the
neighboring agents to control the matching constraints. Messages of
the type $m_{ii'\to f_i}^{(t)}$ and $m_{ii'\to g_{i'}}^{(t)}$ serve
this purpose. They also contain the information about the weights of
the edges (term $\alpha \vw^T\vx$ in the cost function). Similarly,
any two agents that form a square should communicate, so that we can
calculate the term $\beta \vx^T\mS\vx$ in the cost function. This
information is passed by the messages of type $m_{ii'\to
h_{ii'jj'}}^{(t)}$.

\begin{figure}
  \centering
 \includegraphics[width=0.5\columnwidth]{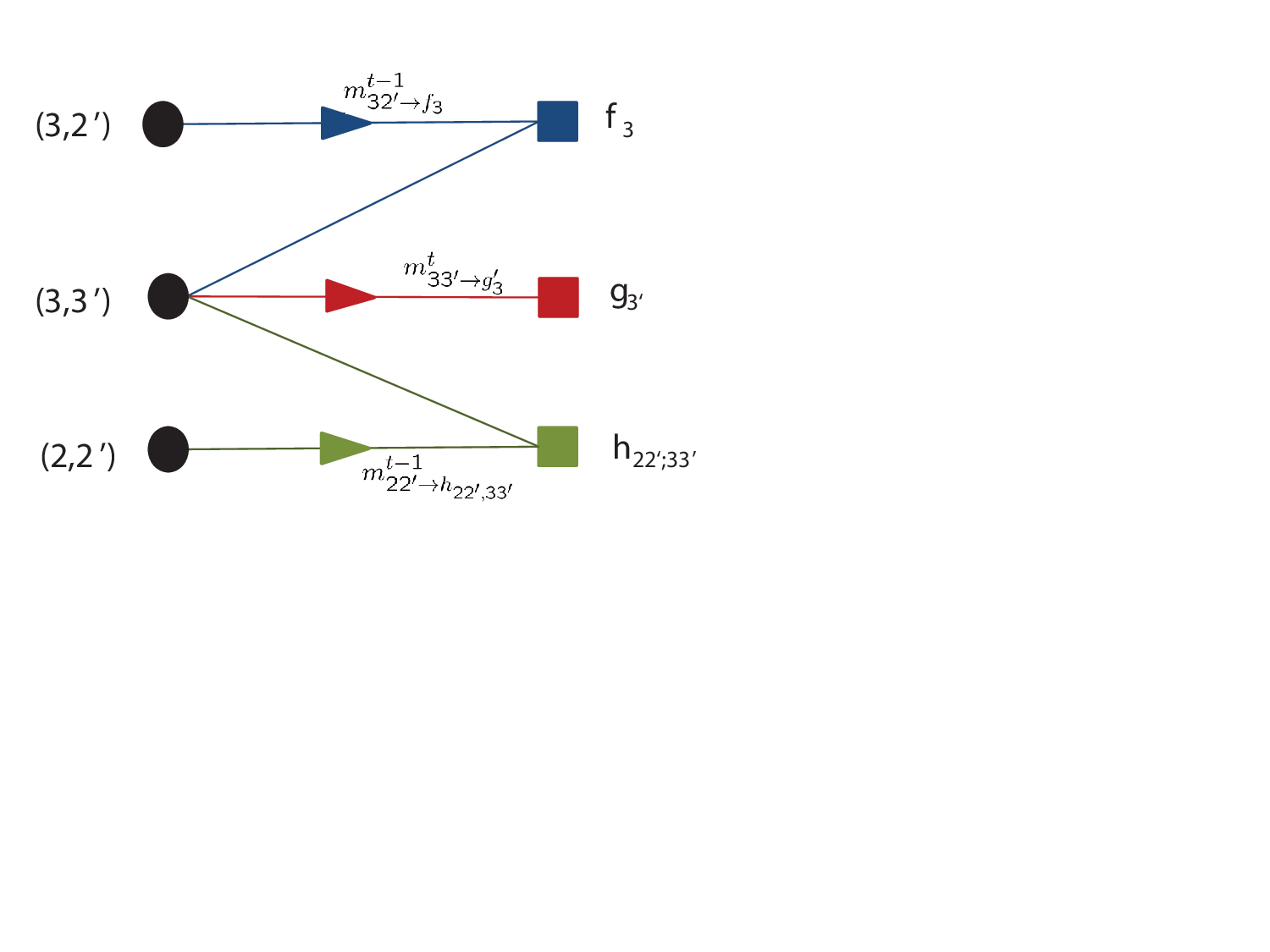}
 \caption{Dependence of $m_{33'\to g_{3'}}^t$ to messages of time $t-1$ for the base example from Figure \ref{fig:netalign-qp-data}.}
   \label{fig:bp-intuition}
\end{figure}

From a slightly different perspective, our algorithm can be seen as
a form of dynamic programming generalized from trees to general
graphs. In fact, it is instructive to consider the special case in
which the factor graph (explained in Section
\ref{subsec:factor-graph}) is indeed a forest. In that case,
removing an edge (or agent) splits the  tree component into two
pieces. This means that the optimization problem \ref{eq:naqp} could
be solved independently on each component. The message of the form
$m_{ii'\to f_i}^{(t)}$ carries the information about the component
that contains $i'$. Figure \ref{fig:bp-intuition} shows this type of
message update. It also contains the information about all squares
that contain $ii'$. Ideally, the message $m_{ii'\to f_i}^{(t)}$
should show the amount of change in the cost function (excluding the
connected component containing $i$) by participation of the edge
$ii'$ in a solution. Similarly, each message of the type $m_{jj'\to
h_{ii'jj'}}^{(t)}$ should be the change in the cost function by
participation of $jj'$
(restricted to the component the edges $jj'$). 

Now we give a rough derivation of equation \eqref{eq:bp1}
using the above discussion. If $ii'$ is
present in the solution, then $\alpha w_{ii'}$ is added
to the cost function. But none of the edges $ki'$ ($k\neq i$)
can now be in the matching. Thus, we should subtract their
maximum contribution
$\left(\max_{k\neq i}\left[m_{ki'\to g_{i'}}^{(t-1)}\right]\right)_+$.
  This explains the first
two terms in the right hand side of equation \eqref{eq:bp1}.
Moreover, we should add the number of squares that will be added by
this edge. For each square $ii'jj'$ if the edge $jj'$ is not present
in the matching, then nothing is added. Otherwise, a $\beta/2$ plus
the term $m_{jj'\to h_{ii'jj'}}^{(t)}$ should be added. This roughly
explains the addition of the third term in \eqref{eq:bp1}. A similar
explanation justifies \eqref{eq:bp2} as well.
%
%
%

%
%
\subsubsection{Convergence of NetAlignMP}\label{sec:BP-Convergence}

We now elaborate on step (4) of \namp.  Ideally, at the end of
iteration $t$, each vertex $i$ selects the edge $ii'$ that sends the
maximum incoming message $m_{ii'\to f_i}^{(t)}$ to it, and we denote
the resulting matching by $M(t)$.  We'd like to terminate the
iteration when $M(t)$ converges.  Unfortunately, picking edges with
this rule does not always produce a matching, and also $M(t)$ may
not converge.  We discuss better approaches to picking a matching
from the messages in Section~\ref{sec:rounding}.  When $M(t)$ does
not converge, it often oscillates between a few states. Therefore,
we could terminate the algorithm when such an oscillation is
observed, and use the current messages to find a matching using the
recipe in Section \ref{sec:rounding}. Another approach for resolving
the oscillation is to use a damping factor $\gamma\in [0,1]$
\cite{MWJ99,BrZ06,frey2007-affinity}. Let $\vec{\mf{n}}(t)$ be the
vector of all messages at time $t$. That is $\vec{\mf{n}}(t)$ is a
fixed ordering of all messages at time $t$. Then the update
equations \eqref{eq:bp1}-\eqref{eq:bp2} can be rewritten as
$\vec{\mf{n}}(t)=F(\vec{\mf{n}}(t-1))$ where $F$ is an operator that
is uniquely defined by equations \eqref{eq:bp1}-\eqref{eq:bp2}. Now
one can consider a new operator $G$ defined by
$G(\vec{\mf{n}}(t))=(1-\gamma^t)\vec{\mf{n}}(t-1)+\gamma^tF(\vec{\mf{n}}(t-1))$
and update the messages using $G$ instead of $F$. The new update
equations will converge for $\gamma<1$.  We make the damping
explicit in the matrix version of this algorithm in
Section~\ref{subsec:bp-matrix-form}.

%

%
%
\SubSection{A matrix formulation} \label{subsec:bp-matrix-form}
We now restate the NetAlignMP algorithm (from
Section~\ref{sec:netalignmp}) using matrix notation.
This helps clarify issues of data organization and computation.
To begin, we again need another bit of notation.
For $\mA \in \mathbb{R}^{m,n}$ and $\vx \in \mathbb{R}^n$,
  define
  \[ \mA \boxdot \vx
     \equiv
     \pmat{
        \max_j a_{1,j} x_j \\
        \max_j a_{2,j} x_j \\
        \vdots \\
        \max_j a_{m,j} x_j
     }. \]
  This operator is just the regular matrix-vector product but with
  the summation $(\mA \vx)_i = \sum_{j} a_{i,j} x_j$ replaced by
  maximization. (This is the matrix-vector product from the
  max-product algebra and is related to the max-plus algebra via
  logarithm/exponential transforms.)
We also need to split
the constraint matrix $\mC$ into $\pmat{\mC_A^T & \mC_B^T}^T$ corresponding
to the matching constraints from graph $A \to B$ and graph $B \to A$,
respectively.

\begin{inlinealgorithm}{NetAlignMP (Matrix-based)}%
 {$\mC = \pmat{\mC_A^T & \mC_B^T}^T$,
  $\mS$, $\vw$, damping parameter $\gamma$, n$_{\textrm{iter}}$, damping\_type}
\begin{lstlisting}[language=Matlab,%
  numbers=left,numberstyle=\scriptsize,
  xleftmargin=5mm,mathescape=true,columns=flexible]
$\vy\itn{0} = 0, \vz\itn{0} = 0, \mS\itn{0} = 0$
for $t$ = 1 to n$_{\textrm{iter}}$
  $\mF = \bound_{0,\frac{\beta}{2}} ({\mS\itn{t-1}}^T + \frac{\beta}{2} \mS)$
  $\vd\itn{t} = \mF \cdot \ve$
  $\vy\itn{t} = \alpha \vw -
               \bound_{0,\infty}[(\mC_A^T \mC_A - \eye) \boxdot \vz\itn{t-1}] +
               \vd\itn{t}$
  $\vz\itn{t} = \alpha \vw -
               \bound_{0,\infty}[(\mC_B^T \mC_B - \eye) \boxdot \vy\itn{t-1}] +
               \vd\itn{t}$
  $\mS\itn{t} = (\mY\itn{t} + \mZ\itn{t} - \alpha \mW - \mD\itn{t})\cdot \mS - \mF$
  if damping_type is 1
    $(\vy\itn{t},\vz\itn{t},\mS\itn{t}) \leftarrow \gamma^t (\vy\itn{t},\vz\itn{t},\mS\itn{t}) + (1-\gamma^t)(\vy\itn{t-1},\vz\itn{t-1},\mS\itn{t-1})$
  else if damping_type is 2
    $\vp = \vy\itn{t-1} + \vz\itn{t-1} - \alpha \vw + \vd\itn{t-1}$
    $(\vy\itn{t},\vz\itn{t},\mS\itn{t}) \leftarrow (\vy\itn{t},\vz\itn{t},\mS\itn{t}) + (1-\gamma^t)(\vp,\vp,\mS\itn{t-1} + {\mS\itn{t-1}}^T - \beta \mS)$
  else if damping_type is 3
    $\vp = \vy\itn{t-1} + \vz\itn{t-1} - \alpha \vw + \vd\itn{t-1}$
    $(\vy\itn{t},\vz\itn{t},\mS\itn{t}) \leftarrow \gamma^t (\vy\itn{t},\vz\itn{t},\mS\itn{t}) + (1-\gamma^t)(\vp,\vp,\mS\itn{t-1} + {\mS\itn{t-1}}^T - \beta \mS)$
  end
  $\vx\itn{t}$ = round_messages($\vy\itn{t},\vz\itn{t},\mS\itn{t}$)
  $\text{obj}\itn{t}$ = objective($\vx\itn{t}$)
end
return $\vx\itn{t}$ with the largest value of $\text{obj}\itn{t}$
\end{lstlisting}
\end{inlinealgorithm}

\subsection{Improved \namp}\label{sec:netalignmp++}

Recall that Klau's algorithm \cite{klau2009-network-alignment} is
obtained by tightening the linear program \ref{eq:nallp}
using combinatorial properties of the problem.
Similarly, we can modify the factor graph
representation of Section \ref{subsec:factor-graph} to
improve the solutions of \namp at
the expense of increasing the running time.  Here is a rough
explanation of this modification. For each variable node $ii'$
add function nodes $d_{ii',j}$ and $d_{ii,j'}$ for all $jj'$
with $ii'\Box jj'$.
These function nodes are defined by:
\begin{align*}
d_{ii',j}([x_{jk'}]_{k':jk'\Box ii'})&=\begin{cases} 1 & \sum_{k': jk'\Box ii'}x_{jk'}\leq 1 \\ 0 & \text{otherwise} \end{cases} \\
d_{ii',j'}([x_{kj'}]_{k: kj'\Box ii'})&=\begin{cases} 1 & \sum_{k: kj'\Box ii'}x_{kj'}\leq 1 \\ 0 & \text{otherwise.} \end{cases}
\end{align*}
After a similarly-lengthy-but-straightforward derivation like
in Appendix~\ref{sec:BP-eqs}, we arrive at the following extension of \namp.

\begin{inlinealgorithm}{\namppp}{%
  $\alpha, \beta$, the set of squares $V_S$, and
  the weighted bipartite graph $\graph{L} = (V_A \cup V_B,E_L)$,
  and a damping parameter $\gamma$.}
\begin{enumerate}
\item At times $t=0,1,\ldots$, each edge $ii'$ sends two messages of the form $m_{ii'\to f_i}^{(t)}$ and $m_{ii'\to g_{i'}}^{(t)}$ and also sends one message of the form $m_{ii'\to h_{ii'jj'}}^{(t)}$ for any square $ii'jj'\in V_S$. {Each square $ii'jj'$ sends a message of the form $m_{ii'jj'\to h_{ii'jj'}}^{(t)}$ and four messages of the type $m_{ii'jj'\to d_{ii',j'}}^{(t)}$ to $d_{ii',j},d_{ii',j'},d_{jj',i}$ and $d_{jj',i'}$.}

\item Messages are initialized by an arbitrary number (let us say $0$).

\item For $t \geq 1$, the messages in iteration $t$ are obtained from the
messages in iteration $t-1$ recursively. In particular for all $ii'\in E_L$
\begin{multline}
\label{eq:bp++1}
m_{ii'\to f_i}^{(t)} = \alpha w_{ii'} -
\left(\max_{k\neq i}\left[m_{ki'\to g_{i'}}^{(t-1)}\right]\right)_+\\
+\sum_{ii'jj'\in V_S} {\left[\left(m_{ii'jj'\to h_{ii'jj'}}^{(t-1)}
+m_{jj'\to h_{ii'jj'}}^{(t-1)}\right)_+ -(m_{jj'\to h_{ii'jj'}}^{(t-1)})_+\right]}.
\end{multline}
The update rule for $m_{ii'\to g_i'}^{(t)}$ is similar to the update rule for $m_{ii'\to f_i}^{(t)}$ and
\begin{multline}
\label{eq:bp++2}
m_{ii'\to h_{ii'jj'}}^{(t)} = \alpha w_{ii'}
+\sum_{\stackrel{kk'\neq jj'}{ii'jj'\in V_S}}
{\left[\left(m_{ii'kk'\to h_{ii'kk'}}^{(t-1)}+m_{kk'\to h_{ii'kk'}}^{(t-1)}\right)_+-(m_{kk'\to h_{ii'kk'}}^{(t-1)})_+\right]}\\
-\left(\max_{k\neq i}\left[m_{ki'\to g_{i'}}^{(t-1)}\right]\right)_+
-\left(\max_{k'\neq i'}\left[m_{ik'\to f_i}^{(t-1)}\right]\right)_+.
\end{multline}
and
\begin{multline}
\label{eq:bp++3}
m_{ii'jj'\to h_{ii'jj'}}^{(t)} = \frac{\beta}{2}
-\left(\max_{k\neq i}\left[m_{ki'jj'\to d_{jj',i'}}^{(t-1)}\right]\right)_+-\left(\max_{k'\neq i'}\left[m_{ik'jj'\to d_{jj',i}}^{(t-1)}\right]\right)_+\\
-\left(\max_{k\neq j}\left[m_{ii'kj'\to d_{ii',j'}}^{(t-1)}\right]\right)_+
-\left(\max_{k'\neq j'}\left[m_{ii'jk'\to d_{ii',j}}^{(t-1)}\right]\right)_+
\end{multline}
and
\begin{multline}
\label{eq:bp++4}
m_{ii'jj'\to d_{ii',j}}^{(t)} = \frac{\beta}{2}
-\left(\max_{k\neq i}\left[m_{ki'jj'\to d_{jj',i'}}^{(t-1)}\right]\right)_+-\left(\max_{k'\neq i'}\left[m_{ik'jj'\to d_{jj',i}}^{(t-1)}\right]\right)_+\\
-\left(\max_{k\neq j}\left[m_{ii'kj'\to d_{ii',j'}}^{(t-1)}\right]\right)_+
+\min\left(m_{ii'\to h_{ii'jj'}}^{(t-1)}+m_{jj'\to h_{ii'jj'}}^{(t-1)},~m_{ii'\to h_{ii'jj'}}^{(t-1)},~m_{jj'\to h_{ii'jj'}}^{(t-1)}\right)
\end{multline}
Equations for $m_{ii'jj'\to d_{ii',j'}}^{(t)}, m_{ii'jj'\to d_{jj',i}}^{(t)}$ and $m_{ii'jj'\to d_{jj',i'}}^{(t)}$ are similar to \eqref{eq:bp++4}.
\item Damp the messages using one of the schemes from NetAlignMP.
\item Round the messages to an integer solution (see possibilities in
Section~\ref{sec:rounding}) and compute the objective function on the
rounded messages
\end{enumerate}
\normalsize {\small \sffamily OUTPUT} the rounded solution with the best objective value.
\end{inlinealgorithm}

%
%

%

%
%

\section{Rounding strategies}\label{sec:rounding}

All algorithms introduced so far rely on formulating
the problem as a mathematical program, with the
integer constraint relaxed. As a consequence, the
computed solution is fractional for most instances.
For IsoRank and Klau's algorithm, the fractional
values are associated with edges in $L$ and for \namp
and \namppp, the values are on both edges and squares.
The last step of each algorithm is to round this
fractional solution to an integral solution, i.e.\ a
matching. There are many ways of rounding, and as always,
the best rounding scheme depends on the
actual problem and the type of relaxation.

The primary type of rounding used is based on
using the fractional solution or the BP messages to
construct a max-weight matching problem. Solving
it produces a solution that then obeys the
matching constraints.  Specifically, we utilize
the function:

\begin{inlinealgorithm}{round\_messages}{messages from $A$ to $B$ $\vy\itn{t}$,
     messages from $B$ to $A$ $\vz\itn{t}$,
     messages on the squares $\mS\itn{t}$}
\begin{lstlisting}[language=Matlab,
  numbers=left,numberstyle=\scriptsize,xleftmargin=5mm,
  mathescape=true,columns=flexible]
$\vx\itn{k}_A = $ bipartite_match(L, $\vy\itn{k}$)
$\text{obj}\itn{k}_A = $ objective($\vx\itn{k}_A$)
$\vx\itn{k}_B = $ bipartite_match(L, $\vz\itn{k}$)
$\text{obj}\itn{k}_B = $ objective($\vx\itn{k}_B$}
return $\vx\itn{k}_{\cdot}$ with the highest value of $\text{obj}\itn{k}_{\cdot}$
\end{lstlisting}
\end{inlinealgorithm}

This function rounds both types of messages and returns the best solution.
Another alternative is to use a greedy matching scheme, where $M$ starts as an empty matching, and
we greedily add edges to $M$ based on the largest values of $\vy\itn{k}_{ii'}$ or $\vz\itn{k}_{ii'}$
such that it stays a matching.
Though computationally more expensive, MWM rounding
yields the best result in most of our experiments.
Therefore, results in Section \ref{sec:synthetic}
and \ref{sec:realdata} are all obtained using MWM
rounding. For the BP algorithm, greedy rounding
using messages on squares -- using $\mS\itn{t}$ above --
yields similar performance as the MWM rounding.

\section{Synthetic experiments} \label{sec:synthetic}

We first compare the belief propagation (BP) algorithm to
existing algorithms on two synthetic matching problems.
The first problem aligns two perturbed grids and the second
aligns two perturbed power-law graphs.

Let $A$ and $B$ be independent realizations of a perturbed
$k \times k$ grid.  The perturbation is a set
of random edges generated with probability $q/d(u,v)^2$
where $d(u,v)$ is the graph distance
between $u$ and $v$.  In these problems, the \emph{ideal alignment}
is known: match each vertex to its image in the other grid.
Now we generate $L$ by matching each grid vertex
to its image and then add additional edges to $L$ with
probability $p$.  This noise globally corrupts the alignment.
We further disturb $L$ by adding
random edges within graph distance $d$ of the end
points of ideal alignment, sampled with probability proportional
to the maximum number of paths.  This step locally corrupts
the alignment.

For the power-law graph test, we construct a reference graph from a
power-law degree sequence with exponent $\theta$ and $n$ vertices
using the algorithm from \citet{Bayati-Saberi-Kim-2007}.
Again, let $A$ and $B$ be
independent realizations of the power-law graph perturbed with the
same noise as the grid above.  Generate $L$ in the same manner, but
without additional distance based edges.

\begin{figure}[t]
 \centering
 \subfigure[Grid graph ($k=20,q=2,d=1$). Function values at left
 and correct matches at right.]{%
   \includegraphics[width=0.45\linewidth]{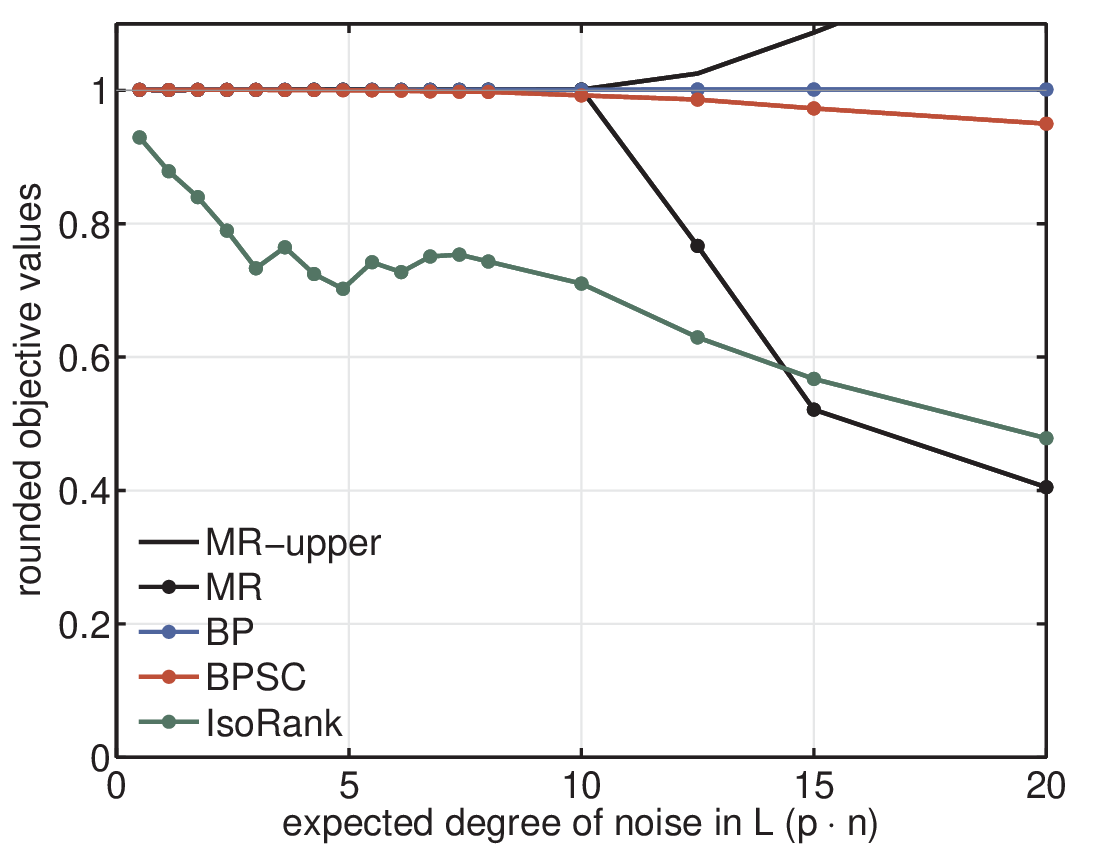}%
   \includegraphics[width=0.45\linewidth]{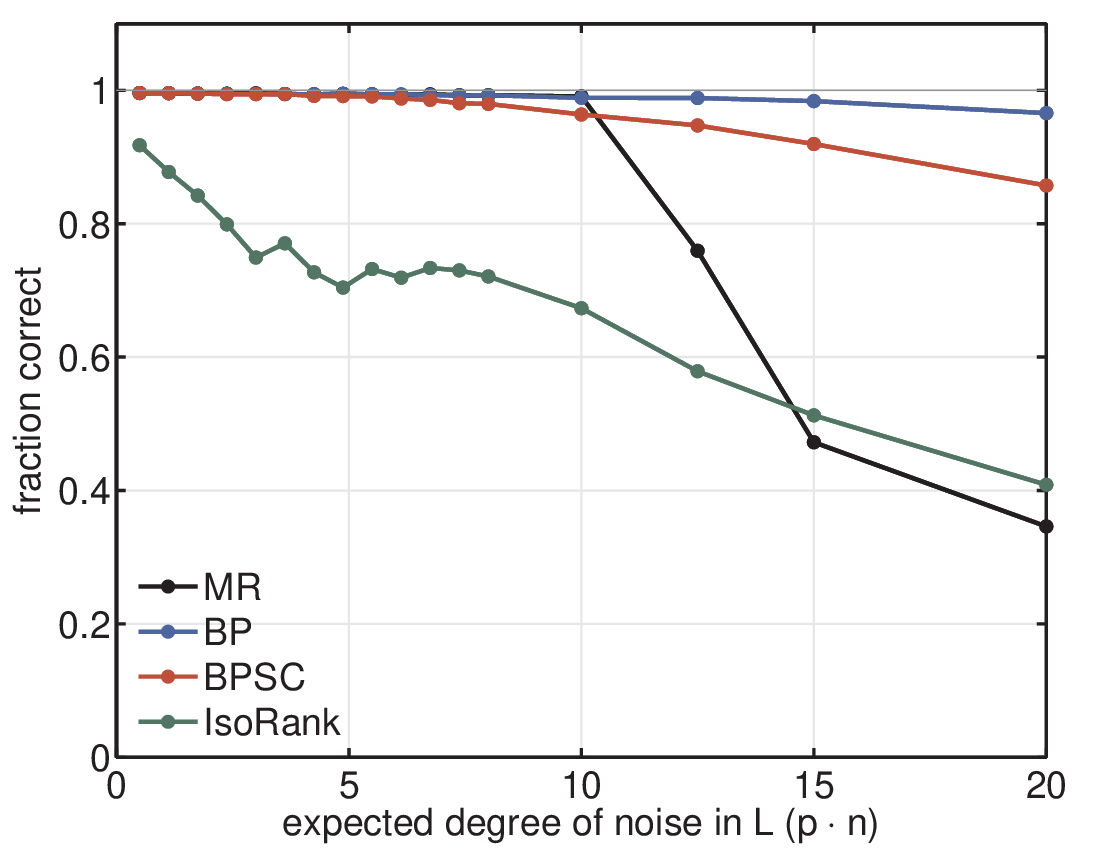}%
 }
  \subfigure[Power-law graphs ($\theta=1.8, n=400,q=1$).
  Function values at left and correct matches at right. function value]{%
   \includegraphics[width=0.45\linewidth]{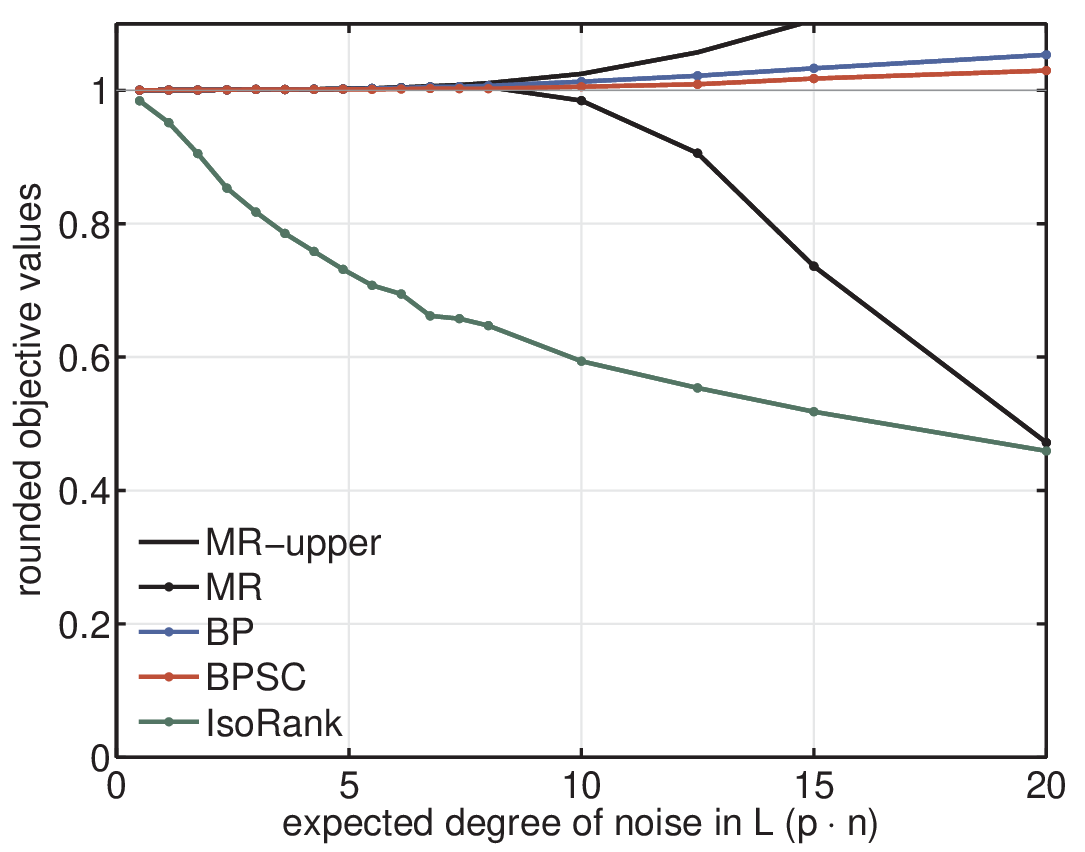}%
   \includegraphics[width=0.45\linewidth]{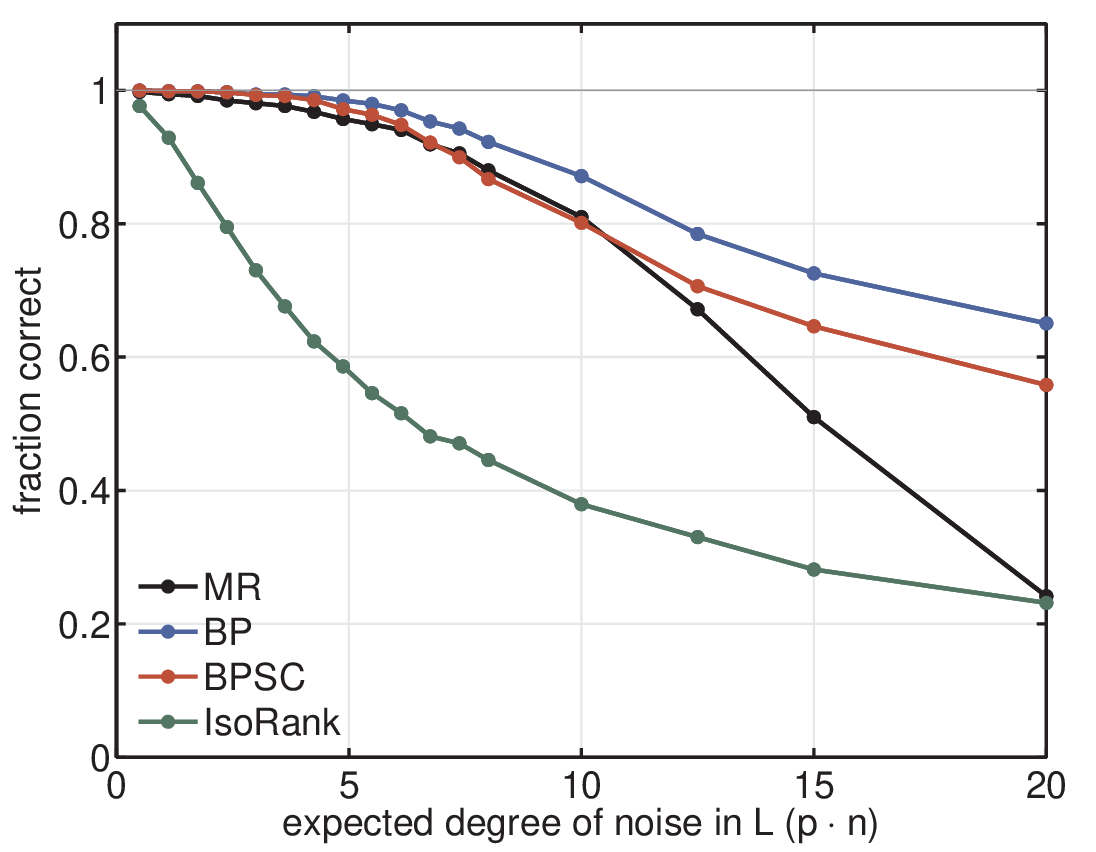}%
 }
 \caption{Upper bounds and correct solutions to synthetic problems
 on grid-graphs (a) and power-law graphs (b).
 The BP and BPSC labels are for the \namp and \namppp algorithms.
 The MR label is for the NetAlignMR algorithm.  The $x$-axis
 corresponds to expected degree that increases with $p$, the fraction of global mismatched
 edges in $L$, which we measure in the expected degree of the noise.
 Once the noise is large, the two message passing approaches
 show the best results.
 Section~\ref{sec:synthetic} for more information.
  }
 \label{fig:synthetic}
\end{figure}

In our results, we compare all outputs to the reference matching between
the graphs $A$ and $B$.  Figure~\ref{fig:synthetic} shows
the average fraction of the reference matching
obtained by each algorithm over 48 trials.  The objective
function is pure overlap and the dark lines in the figure show the
ratio of the algorithm's overlap to the overlap of the reference solution.
Each algorithm should be computing a good objective, and thus larger values
are better.  Indeed, the reference solution may not be the best solution
when $L$ is highly corrupted with a large expected degree.
When this happens with the power-law graphs, we observe
that the BP algorithm finds a matching with a higher overlap and thus
the fraction is larger than 1.
Similarly, the light lines show the fraction of correct matches
from the the algorithms.  These values track the objective values
showing that the network alignment objective is a good surrogate for
the number of correct matches objective.

When the amount of random noise in $L$ exceeds an expected degree of $10$
for the grid graphs and 8 for the power-law graphs,
many of the algorithms are no longer able to obtain good solutions.
In this regime, the MP and \mppp algorithms performs better than
the MR algorithm.

We used the MP and \mppp algorithms with $\alpha=1,\beta=2$, the IsoRank algorithm
with $\gamma = 0.95$, and the MR algorithm with $\alpha=0,\beta=1$
for these experiments. These parameters are natural for the various
algorithms.  For example, MR requires $\alpha=0, \beta=1$ to produce
an upper-bound on overlap. Below, we study the behavior of the algorithms
for a wider variety of parameters.

%

\Section{Real datasets} \label{sec:realdata}

While we saw that the BP algorithm performed well on noisy synthetic problems
in the previous section, in this section we investigate alignment problems
from bioinformatics and ontology matching.  For each algorithm, we
explore a range of choices for all of the parameter values and summarize
the results from the \emph{best} choice in Table~\ref{tab:summary}.
Note that Klau's algorithm uses two parameters $\gamma$ and
$st$ to control the subgradient method.

\begin{table}
\tbl{Properties of the real-world test problems.}{
\footnotesize
\begin{tabularx}{\linewidth}{lXXXXX}
\toprule
Problem & $|V_A|$ & $|E_A|$ & $|V_B|$ & $|E_B|$ & $|E_L|$
\\ \midrule
    dmela-scere    & 9459 & 25636 & 5696 & 31261 & 34582
\\  Mus M.-Homo S. & 3247 & 2793 & 9695 & 32890 & 15810
\\  lcsh2wiki-small & 1919 & 1565 & 2000 & 3904 & 16952
\\  lcsh2wiki-full  & 297266 & 248230 & 205948 & 382353 & 4971629
\\ \bottomrule
\end{tabularx}
}
\end{table}

\SubSection{Bioinformatics} \label{subsec:bioinfo}

The alignment of protein-protein interaction (PPI) networks of
different species is an important problem in
bioinformatics~\cite{singh2007-matching-topology}. We
consider aligning the PPI network of Drosophila melanogaster (fly)
and Saccharomyces cerevisiae (yeast), and Homo sapiens (human) and
Mus musculus (mouse). These PPI networks are available in several
open databases and they are used in \cite{singh2008-isorank-multi}
and \cite{klau2009-network-alignment}, respectively. For each
problem, we utilize the value of $\vw$ from the
original publication.  While the
results of the experiment are rich in biological information,
we focus solely on the optimization problem.

\begin{figure}
  \centering
 \includegraphics[width=0.9\linewidth]{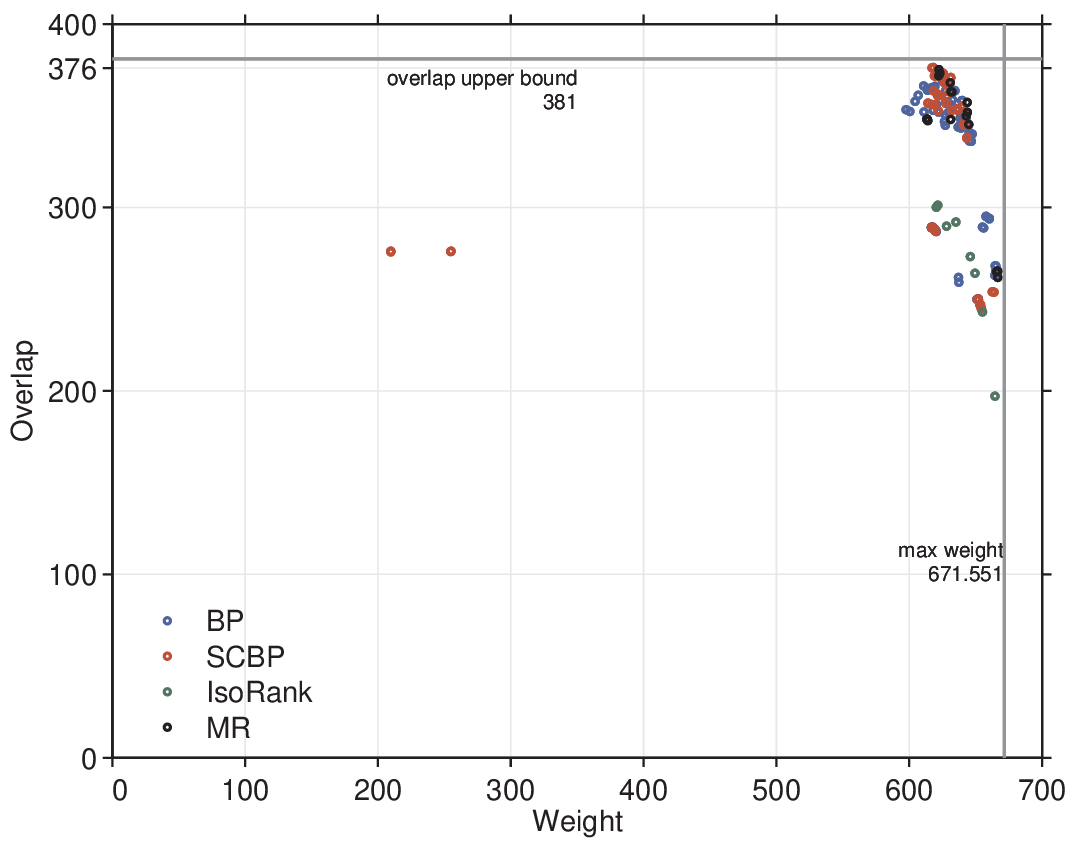}
 \hfill
 \includegraphics[width=0.9\linewidth]{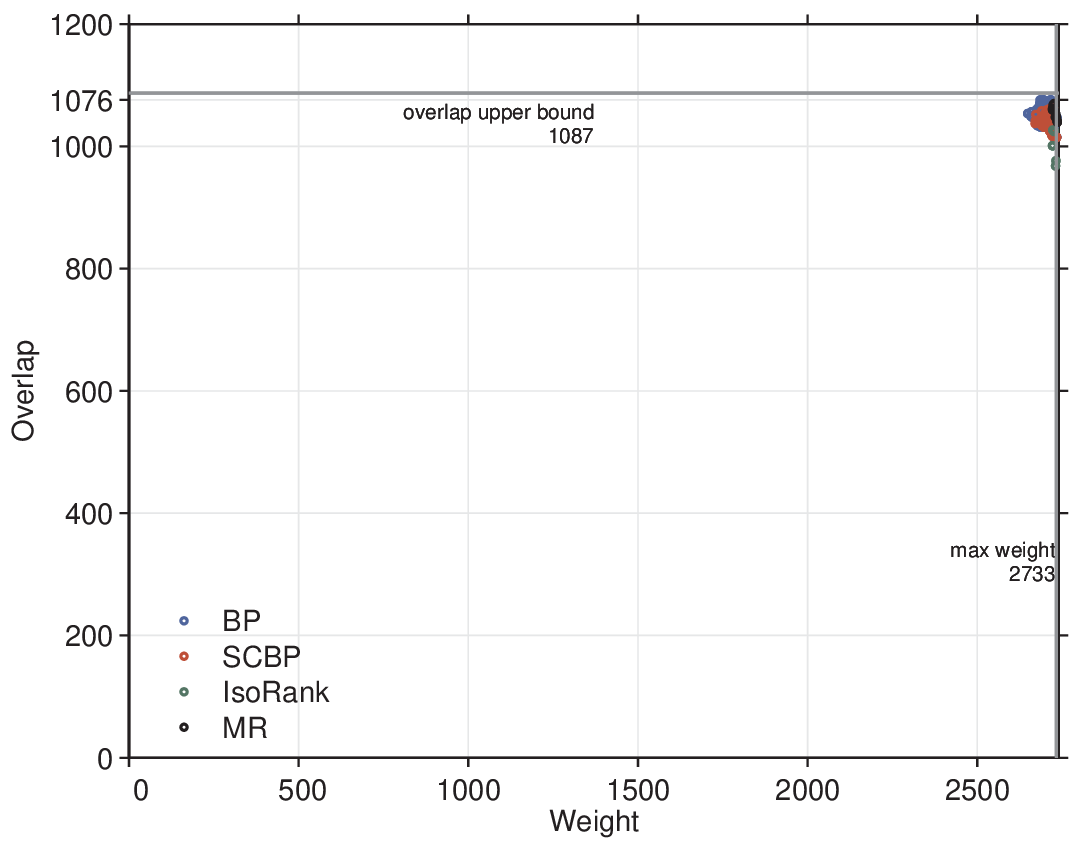}
 \caption{Results of the three algorithms SpaIsoRank (IsoRank), \namp (BP),
   \namppp (SCBP) and NetAlignMR (MR) on the Mus. M.-Homo S. alignment
       (top) and dmela-scere alignment (bottom).
   }
   \label{fig:bioinfo}
\end{figure}

Figure~\ref{fig:bioinfo} shows the performance of the four
algorithms  -- \namp, \namppp, NetAlignMR, SpaIsoRank --
on these two alignments.  For each algorithm, we perform a parameter
sweep over the following parameters
\begin{center}
 \begin{tabular}{lll}
  IsoRank & Damping & $\gamma \in \{ 0.3, 0.5, 0.85, 0.95 \}$ \\
          & Rounding & type $ \in \{ 1, 2 \}$ \\
  \midrule
  \namp & Objective & $(\alpha,\beta) \in \{ (10,1), (2,1), (1,1), (1,2), (1,10) \}$ \\
        & Damping & $\gamma \in \{ 0.9, 0.99, 0.995, 0.999 \}$ \\
        &         & type $ \in \{ 2, 3 \}$ \\
  \midrule
  \namppp & \multicolumn{2}{l}{\emph{same as } \namp}\\
  \midrule
  NetAlignMR & Objective & $(\alpha,\beta) \in \{ (10,1), (2,1), (1,1), (1,2), (1,10) \}$ \\
             & Damping & $\gamma \in \{ 0.1, 0.4 \}$ \\
             &          & mstep $ \in \{ 5, 25, 50 \}$ \\
 \end{tabular}
\end{center}
We run IsoRank until convergence, and run the other approaches
for a total of 500 iterations.
On these instances, we
record the best iterate ever generated and plot the overlap
and weight of the alignments in the figure.

In both problems \namp, \namppp and NetAlignMR manage to obtain
near-optimal solutions. In terms of the largest overlap,
our \namp does the best on the Mus M.-Homo S. alignment,
whereas \namppp does the best of the dmele-scere alignment.
(See Table~\ref{tab:summary} for the parameters that
produced the best overlap.)

\SubSection{Ontology} \label{sec:lcsh2wiki}

\begin{figure}
 \centering
 \includegraphics[width=0.9\linewidth]{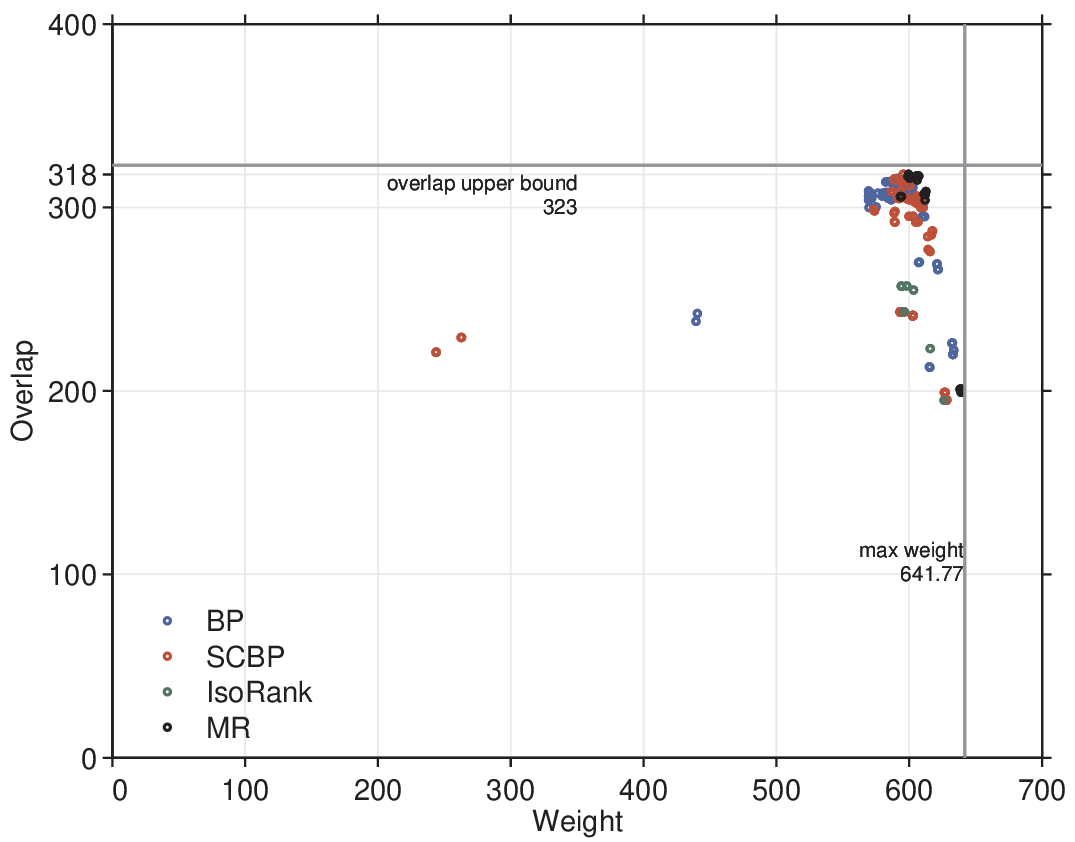}
 \hfill
 \includegraphics[width=0.9\linewidth]{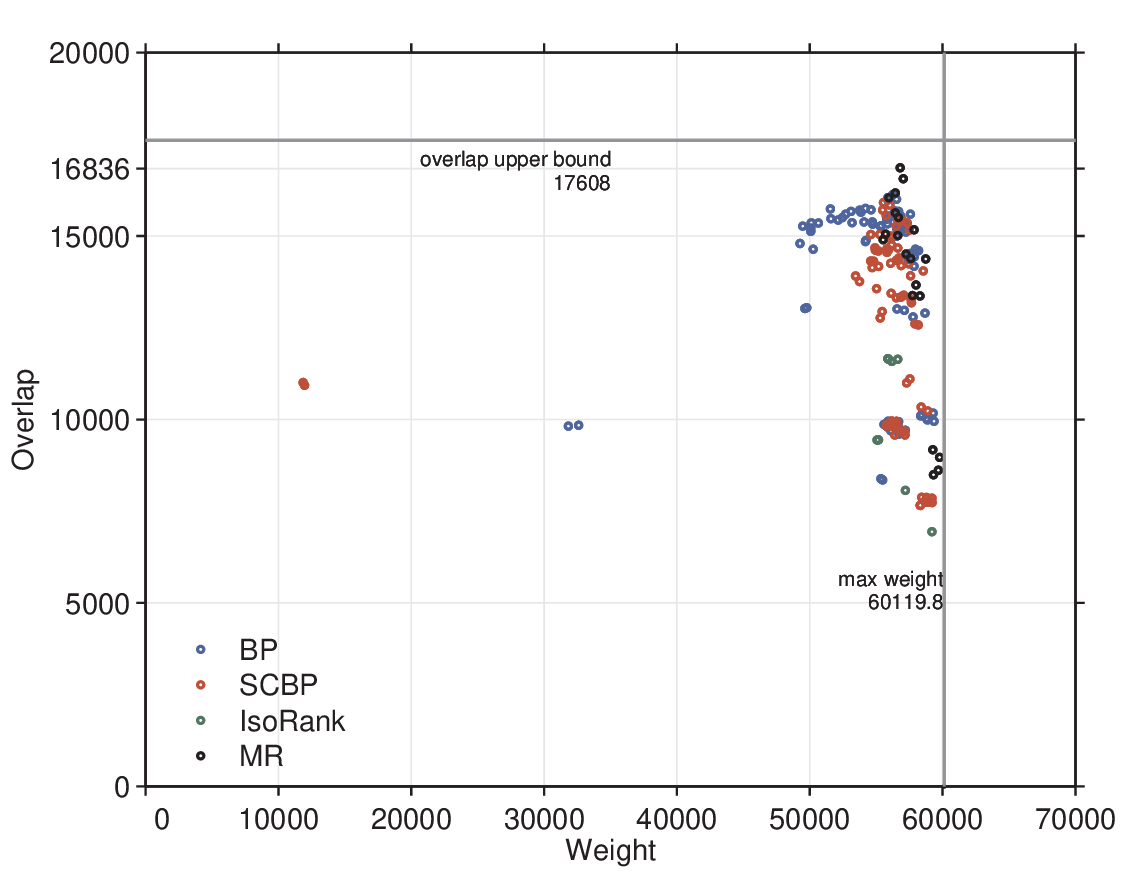}
 \caption{Results of the four algorithms SpaIsoRank (IsoRank), \namp (BP),
   \namppp (SCBP), NetAlignMR (MR) on lcsh2wiki-small (top) and lcsh2wiki (bottom).}
  \label{fig:lcsh2wiki}
\end{figure}

Our original motivation for investigating network alignment is
aligning the network of subject headings from the Library of
Congress with the categories from
Wikipedia~\cite{wikipedia2007-april-data}. Each node in these
networks has a small text label, and we use a Lucene search index
\cite{lucene2007} to quickly find potential matches between nodes
based on the text.
To score the matches, we use the SoftTF-IDF
scoring routine~\cite{cohen2003-matching-names}.  These scores
become the weights in $\vw$.  Our \emph{real}
problem is to match the entire graphs.  From this problem we extract
a small instance that should capture the most important nodes in the
problem. (Node importance is either reference count (subject
headings) or PageRank (Wikipedia categories).) The results are shown
in Figure~\ref{fig:lcsh2wiki}.

We repeated the parameter sweep from the previous section on these
two problems as well.  The best algorithm on these two problems
is NetAlignMR, with \namp and \namppp alternating for second place.
In lcsh2wiki-small, the upper bound computed by NetAlignMR is 323.
\namp achieves a lower bound of 318 and NetAlignMR achieves
321. In lcsh2wiki, we compute an upper bound of 17608 using
a linear programming solver on \ref{eq:natlp} with the full symmetry
constraints instead of the Lagrange multipliers.
Though not shown in the
figure, \namp obtains a lower bound of 16204 with $\gamma =
0.9995$, $\alpha = 0.2$ and $\beta = 1$.

In all our real datasets $L$ is quite sparse, making NetAlignMR more
favorable. Still, \namp is closely following and has an
advantage on running time -- see the summary in Table~\ref{tab:summary}
for information about runtime.

\begin{table}
\tbl{The parameters used to produce the results with
the highest overlap from Figures~\ref{fig:bioinfo}~and~\ref{fig:lcsh2wiki}.
We abbreviated lcsh2wiki-small as lcsh-small.  The overlap score
shows the highest overlap produced by that method on the problem
and the percentage of the best upper-bound on the solution objective.
All times are reported in seconds, and the Sol. Time column indicates
the time taken to compute the best solution whereas the Total Time
column indicates the total time for all iterations of the method.
(Recall that all methods return the iterate with the best solution,
which may not be the final iterate.)
We organized the table to indicate the most successful parameter choices.
}
{
\footnotesize
\begin{tabularx}{\columnwidth}{@{}llr@{\;\;}lXXl@{}}
\toprule
Alg. & Data & \multicolumn{2}{r}{Overlap} & Sol. Time & Total Time & Parameters\\
\midrule
MWM &       musm-homo &     393 &  36.2\% &     --- &     --- & \\
    &     dmela-scere &     135 &  35.4\% &     --- &     --- & \\
      & lcsh-small &     119 &  36.8\% &     --- &     --- & \\
      &       lcsh2wiki &    2346 &  13.3\% &     --- &     --- & \\
\midrule
Iso &       musm-homo &    1027 &  94.5\% &     0.0 &     0.4 & $\gamma= 0.50$; r=2 \\
      &     dmela-scere &     301 &  79.0\% &     3.7 &    10.7 &
$\gamma=0.95$; r=2 \\
      & lcsh-small &     257 &  79.6\% &     0.0 &     0.7 &
$\gamma=0.50$; r=2 \\
      &       lcsh2wiki &   11732 &  66.6\% &    11.7 &   587.3 &
$\gamma=0.95$; r=2 \\
\midrule
MP &       musm-homo &    1076 &  99.0\% &     2.6 &    13.2 &
$\alpha=2;\beta=1;\gamma=0.995;d=3$ \\
      &     dmela-scere &     369 &  96.9\% &    26.7 &    34.9 &
$\alpha=1; \beta=2; \gamma=0.999;d= 3$ \\
      & lcsh-small &     316 &  97.8\% &     7.6 &    12.6 &
$\alpha=1; \beta=  1; \gamma= 0.999; d= 3$ \\
      &       lcsh2wiki &   15974 &  90.7\% &  3795.3 &  4198.4 &
$\alpha =1; \beta=  2; \gamma= 0.999; d= 2$ \\
\midrule
\mppp &       musm-homo &    1062 &  97.7\% &    14.4 &    17.3 &
$\alpha=  1; \beta=  1; \gamma= 0.999; d= 3$ \\
      &     dmela-scere &     376 &  98.7\% &    28.7 &    33.3 &
$\alpha= 1; \beta= 10; \gamma= 0.999; d= 3$ \\
      & lcsh-small &     318 &  98.5\% &    11.8 &    15.2 &
$\alpha= 1; \beta=  2; \gamma =0.999; d=3$ \\
      &       lcsh2wiki &   15771 &  89.6\% &  4103.8 &  4990.2 &
$\alpha=  1; \beta=  1; \gamma= 0.999; d= 3$ \\
\midrule
MR&       musm-homo &    1070 &  98.4\% &    12.5 &    12.6 &
$\alpha = 1; \beta =10; \gamma= 0.400; st= 5$ \\
      &     dmela-scere &     375 &  98.4\% &    22.7 &    79.4 &
$\alpha=  1; \beta = 2; \gamma= 0.400; st= 5$ \\
      & lcsh-small &     318 &  98.5\% &     4.1 &    16.8 &
$\alpha = 1; \beta = 2; \gamma= 0.400; st= 5$ \\
      &       lcsh2wiki &   16836 &  95.6\% &  4878.2 &  4988.0 &
$\alpha=  1; \beta = 2; \gamma =0.400; st= 5$ \\
\bottomrule
\end{tabularx}
}
\label{tab:summary}
\end{table}

\SubSection{Multi-lingual ontologies} \label{lcsh2rameau}
For a final test, we evaluate automatically aligning two large networks
where a correct alignment exists.  The networks
are the Library of Congress Subject Headings and its French
analogue, Rameau.  Both are similar ontologies and we expect
a non-trivial alignment between the networks.  The correct
alignment between the networks is available from
\url{http://www.cs.vu.nl/STITCH/rameau/dump/}.
It contains $57,645$ matches between the $154,974$
nodes of Rameau and the $342,684$ nodes of LCSH.
(This experiment used a newer version of LCSH than the previous
experiments, which is why the number of nodes changed).

To build the set of potential matches, we translate the
French subject headings to English using Google Translate
(translate.google.com),
and translate the English headings to French also using Google
Translate.  Then, we use Lucene to compute a pairwise
match between the strings and keep the top 25 matches.
This produces up to 100 potential matches per node,
25 from LCSH $\to$ Rameau in English, 25 from Rameau $\to$
LCSH in English, and another 50 for the same sets in
French.
The weights are computed in the same way as in the
previous section.
In total, we had $20,883,500$ possible
edges between the graphs.  Of these, only $42,215$
of the correct matches appeared.  The overlap induced
by the correct set of matches is $39,749$.

The results and running time from our algorithms are presented
in Table~\ref{tab:lcsh2rameau}.  In summary, NetAlignMR computes
the best results in terms of the optimization objective,
but it also takes the most time.  \namppp is
the runner-up and fills the gap in results and run-time
between NetAlignMP and NetAlignMR.  With respect
to recall and precision, NetAlignMR has the highest
recall (57.6\%) with good precision (27.0\%), but
\namp and \namppp always have slightly higher precision.
Note that we
performed no specific tuning to
account for the differences in French and English.

In the table, we also showed the number of triangles
overlapped by a matching.
This number appears to be indicative of the true
matching performance. We
believe these results demonstrate that including
overlapped triangles into the objective may
improve matching algorithms.  We plan to investigate
these ideas in the future.

\begin{table}
\tbl{The alignment results for LCSH and Rameau.  The first
set of results shows the statistics of the known alignment
and the results from the max-weight matching algorithm.
Next we show results from our algorithms for three objective parameters.
The columns are: objective parameters, algorithms, matching weight,
matching edge overlap, time, total correct, recall,  precision,
and matching triangle overlap.
}{
\footnotesize
 \begin{tabularx}{\linewidth}{llXXXXXXX}
 \toprule
 Obj. & Alg. & Weight & Overlap &  Time \rlap{(s)} & Correct & Rec. & Prec. & Triangles \\
  \midrule
& Sol. & 36332.42& 39847 & ---&57645 & 100\% & 100\% & 2073 \\
&  MWM & 93279.0 &   16990 &   29.6 &   29098 &  50.5\% &  23.3\% &     350 \\
\midrule
$\alpha=1, \beta=1$ &   MP & 84622.0 &   46400 & 23522.0 &   32585 &  56.5\% &  27.6\% &    1515 \\
                & MP++ & 85810.1 &   46942 & 27115.6 &   32857 &  57.0\% &  27.4\% &    1548 \\
                &   MR & 87588.6 &   48367 & 33366.9 &   33225 &  57.6\% &  27.0\% &    1617 \\
\midrule
$\alpha=1, \beta=2$ &   MP & 81752.6 &   46569 & 23427.1 &   31724 &  55.0\% &  27.6\% &    1483 \\
                & MP++ & 84615.7 &   46656 & 26673.1 &   31952 &  55.4\% &  26.7\% &    1531 \\
                &   MR & 85438.4 &   48934 & 56961.6 &   32303 &  56.0\% &  26.3\% &    1604 \\
\midrule
$\alpha=0, \beta=1$ &   MP & 60617.9 &   45247 & 14284.8 &   24794 &  43.0\% &  23.2\% &    1467 \\
                & MP++ & 60502.8 &   41592 & 13979.5 &   24498 &  42.5\% &  23.0\% &    1484 \\
                &   MR & 65994.2 &   46163 & 10384.4 &   25455 &  44.2\% &  21.5\% &    1602 \\

\bottomrule
 \end{tabularx}
}
\label{tab:lcsh2rameau}

\end{table}

\Section{Discussion}

Let us recap.  Network alignment is an important tool in
a variety of applications including systems biology,
computer vision and ontology matching. It is especially useful
for comparing large datasets with inherent and
related graph structures. Here, we explored
matching protein-protein interaction networks
and ontologies.  In the future, we envision
applications of these techniques in
mapping large social network structure.

Of course, finding the best alignment between two
networks is NP-hard.  Thus far, we are limited
to attacking the problem heuristically as
there is no known approximation algorithm.
Many different heuristics for the problem fit
nicely within our quadratic programming framework
for the problem.
We studied several existing algorithms this framework
and compared their performance on both synthetic
and real data.

We find that the NetAlignMR from \citet{klau2009-network-alignment}
produces the best results when a sparse set
of potential matches between two graphs exist.
Our two new message-passing
algorithms, \namp and \namppp, were designed
based on belief propagation ideas for
solving the integer optimization problem
directly.  They are mildly faster than NetAlignMR
(roughly $1.3\%$ in our experiments) and their results
nearly tie with NetAlignMR.  Additionally,
our algorithms produce better solutions
when the set of potential matches is dense.

There are a number of avenues for future work
we plan to investigate. First, because our algorithms
use message passing, they
should allow simple parallel implementations,
including on MapReduce style architectures.
Second, in each of the real data sets we used,
the nodes of the two graphs had an informative
label, which helped us to apply preprocessing
to produce a sparse graph of potential
matches between the two graphs.  All of the
previously discussed algorithms utilize
this fact, except for IsoRank.
We also plan to investigate aligning graphs
without these initial ``hints.''

\begin{acks}
We thank  Margot Gerritsen for helping with an initial version of
this manuscript.  And we extend our heartfelt thanks to the people
at the Library of Congress
for funding this work, and the computational approaches to digital stewardship
group at Stanford.  In particular, we'd like to thank Laura Campbell,
Barbara Tillet, and Ed Summers for their own contributions.

Thanks to Jure Leskovec for discussing the problem with us.
Also, we thank Nathan Sakunkoo for implementing the Soft TF-IDF
based scoring algorithm.
\end{acks}

\bibliographystyle{acmtrans}
\bibliography{all-bibliography}

\received{February 2007}{March 2009}{June 2009}

\newpage

\appendix

\Section{Derivation of \namp equations}\label{sec:BP-eqs}
{\footnotesize
The belief propagation algorithm (and its max-product version) is an iterative procedure for passing message along the edges of a factor graph \cite{Pea88}.
We use the notation $t=0,1,\ldots$ to denote the messages after $t$ message passing steps.  The BP algorithm specifies the messages to pass for a general factor graph.  For our factor-graph representation we obtain two types of \emph{real-valued} BP messages.
We denote these two types by $\nu$ and $\la$ respectively.
\begin{enumerate}
\item Messages from variable nodes to function nodes. Each variable node $ii'$ sends the following messages
\begin{align}
\nu_{ii'\to f_i}^{(t+1)}(x_{ii'})&=\la_{g_{i'}\to ii'}^{(t)}(x_{ii'})\prod_{ii'jj'}\la_{h_{ii'jj'}\to ii'}^{(t)}(x_{ii'}),\label{eq:m(ij->fi)}\\
\nu_{ii'\to g_{i'}}^{(t+1)}(x_{ii'})&=\la_{f_i\to ii'}^{(t)}(x_{ii'}) \prod_{ii'jj'}\la_{h_{ii'jj'}\to ii'}^{(t)}(x_{ii'}),\label{eq:m(ij->gj)}\\
\nu_{ii'\to h_{ii'jj'}}^{(t+1)}(x_{ii'})&=\la_{f_i\to ii'}^{(t)}(x_{ii'})\la_{g_{i'}\to ii'}^{(t)}(x_{ii'}) \prod_{ii'kk'\neq ii'jj'}\la_{h_{ii'kk'}\to ii'}^{(t)}(x_{ii'}). \label{eq:m(ij->hijrs)}
\end{align}
Note that each variable node $ii'jj'$ has only one neighbor. Hence, its message is always defined by $\nu_{ii'jj'\to h_{ii'jj'}}^{(t+1)}(x_{ii'jj'})=1$.

\item The messages from function nodes to variable nodes are:
\begin{align}
\la_{f_i\to ii'}^{(t)}(x_{ii'})&=\max_{\vx_{\pa f_i\setminus \{ii'\}}}\left\{e^{[\alpha\sum_{j'} w_{ij'}x_{ij'}]}f_i(\vx_{\pa f_i})
\prod_{j'\neq i'}\nu_{ij'\to f_i}^{(t)}(x_{ij'})\right\},\no \\
\la_{g_{i'}\to ii'}^{(t)}(x_{ii'})&=\max_{\vx_{\pa g_{i'}\setminus \{ii'\}}}\left\{e^{[\alpha\sum_{j} w_{ji'}x_{ji'}]}g_{i'}(\vx_{\pa g_{i'}})
\prod_{j\neq i}\nu_{ji'\to g_{i'}}^{(t)}(x_{ji'})\right\}, \no\\
\la_{h_{ii'jj'}\to ii'}^{(t)}(x_{ii'}) &= \max_{x_{jj'},x_{ii'jj'}}\left\{ e^{\frac{\beta}{2} x_{ii'jj'}} h_{ii'jj'}(\vx_{\pa  h_{ii'jj'}}) \nu_{jj'\to h_{ii'jj'}}^{(t)}(x_{jj'})\right\},\no\\
\la_{h_{ii'jj'}\to jj'}^{(t)}(x_{jj'}) &= \max_{x_{ii'},x_{ii'jj'}}\left\{ e^{\frac{\beta}{2} x_{ii'jj'}} h_{ii'jj'}(\vx_{\pa  h_{ii'jj'}}) \nu_{ii'\to h_{ii'jj'}}^{(t)}(x_{ii'})\right\},\no\\
\la_{h_{ii'jj'}\to ii'jj'}^{(t)}(x_{ii'jj'}) &= \max_{x_{ii'},x_{jj'}}\left\{ e^{\frac{\beta}{2} x_{ii'jj'}} h_{ii'jj'}(\vx_{\pa  h_{ii'jj'}}) \nu_{ii'\to h_{ii'jj'}}^{(t)}(x_{ii'}) \nu_{jj'\to h_{ii'jj'}}^{(t)}(x_{jj'})\right\}.\no\\
&\label{eq:m(var->fun)}
\end{align}
\end{enumerate}
At the end of each iteration $t$, each variable node $x_{ii'}$ ($x_{ii'jj'}$) is assigned a binary value as follows:
\begin{align*}
x_{ii'}^{(t)}&=\arg\max_{x_{ii'}}\left\{ \la_{f_i\to ii'}^{(t)}(x_{ii'})\la_{g_{i'}\to ii'}^{(t)}(x_{ii'})\prod_{ii'jj'}\la_{h_{ii'jj'}\to ii'}^{(t)}(x_{ii'})\right\},\\
x_{ii'jj'}^{(t)}&=\arg\max_{x_{ii'jj'}}\left\{ \prod_{ii'jj'}\la_{h_{ii'jj'}\to ii'jj'}^{(t)}(x_{ii'jj'})\right\}.
\end{align*}
In many applications as $t\to \infty$ the assigned values $x_{ii'}^{(t)}, x_{ii'jj'}^{(t)}$ converge to good approximate solutions.

It is possible to simplify the equations above by eliminating
redundancies -- for example, we already mentioned that the
message $\nu_{ii'jj'\to h_{ii'jj'}}^{(t+1)}(x_{ii'jj'})=1$ always.
We now simplify the above set of equations. Since the variables $x_{ij}$ and $x_{ijrs}$ are binary valued, we compress the messages by sending just the log-likelihood values $m_{ij\to f_i}^{(t)}=\log \big(\nu_{ij\to f_i}^{(t)}(1)/\nu_{ij\to f_i}^{(t)}(0)\big)$. Similarly we define messages $m_{ij\to g_j}^{(t)}$, $m_{ij\to h_{ijrs}}^{(t)}$, and $m_{ijrs\to h_{ijrs}}^{(t)}$.

    Next, we will carry out these calculations for $m_{ij\to f_i}^{(t)}$.
\begin{align*}
m_{ii'\to f_i}^{(t+1)}&=\log\left(\f{\la_{g_{i'}\to ii'}^{(t)}(1)}{\la_{g_{i'}\to ii'}^{(t)}(0)}\right)+\sum_{ii'jj'}\log\left(\f{\la_{h_{ii'jj'}\to ii'}^{(t)}(1)}{\la_{h_{ii'jj'}\to ii'}^{(t)}(0)}\right).
\end{align*}
Each $\log$ term here can be simplified because $\log$ is a monotone function, and hence, it commutes with $\max$.  For example,
\begin{align*}
\log(\la_{g_{i'}\to ii'}^{(t)}(1))&=\max_{\stackrel{\vx_{\pa g_{i'}\setminus \{ii'\}}}{x_{ii'}=1}}\left\{\alpha w_{ii'}+\alpha\sum_{j\neq i} w_{ji'}x_{ji'}+\log g_{i'}(\vx_{\pa g_{i'}})+\sum_{j\neq i}\log\nu_{ji'\to g_{i'}}^{(t)}(x_{ji'})\right\}\no\\
&=\alpha w_{ii'}+\sum_{j\neq i}\log\nu_{ji'\to g_{i'}}^{(t)}(0)
\end{align*}
where the last equality uses the matching constraint imposed by $g_{i'}$. Similarly,
\begin{align*}
\log(\la_{g_{i'}\to ii'}^{(t)}(0))&=\max\left\{\sum_{j\neq i}\nu_{ji'\to g_{i'}}^{(t)}(0),\max_{k\neq i}\left[\alpha w_{ki'}+\sum_{j\neq i}\log\nu_{ji'\to g_{i'}}^{(t)}(0)+\log(\f{\nu_{ki'\to g_{i'}}^{(t)}(1)}{\nu_{ki'\to g_{i'}}^{(t)}(0)})\right]\right\}.
\end{align*}
Therefore, we have
\begin{align*}
\log\left(\f{\la_{g_{i'}\to ii'}^{(t)}(1)}{\la_{g_{i'}\to ii'}^{(t)}(0)}\right)&=\alpha w_{ii'}-\left\{\max_{k\neq i}(\alpha w_{ki'}+m_{ki'\to g_{i'}}^{(t)})\right\}_+
\end{align*}
where $(a)_+$ means $\max(a,0)$.  Similar calculations for $\la_{h_{ii'jj'}\to ii'}^{(t)}$ yield
\begin{multline*}
\log\left(\f{\la_{h_{ii'jj'}\to ii'}^{(t)}(1)}{\la_{h_{ii'jj'}\to ii'}^{(t)}(0)}\right)=\max\left(\frac{\beta}{2}+\log\nu_{jj'\to h_{ii'jj'}}^{(t)}(1),\log\nu_{jj'\to h_{ii'jj'}}^{(t)}(0)\right)\\
-\max\left(\log\nu_{jj'\to h_{ii'jj'}}^{(t)}(1),\log\nu_{jj'\to h_{ii'jj'}}^{(t)}(0)\right)\\
=(\frac{\beta}{2}+m_{jj'\to h_{ii'jj'}}^{(t)})_+-(m_{jj'\to h_{ii'jj'}}^{(t)})_+.
\end{multline*}
Summarizing, we obtain
\begin{align*}
m_{ii'\to f_i}^{(t+1)}&=\alpha w_{ii'}-\left\{\max_{k\neq i}(\alpha w_{ki'}+m_{ki'\to g_{i'}}^{(t)})\right\}_++\sum_{ii'jj'}\left((\frac{\beta}{2}+m_{jj'\to h_{ii'jj'}}^{t})_+-(m_{jj'\to h_{ii'jj'}}^{(t)})_+\right).
\end{align*}
By symmetry we obtain
\begin{align*}
m_{ii'\to g_{i'}}^{(t+1)}&=\alpha w_{ii'}-\left\{\max_{k'\neq i'}(\alpha w_{ik'}+m_{ik'\to f_{i}}^{(t)})\right\}_++\sum_{ii'jj'}\left((\frac{\beta}{2}+m_{jj'\to h_{ii'jj'}}^{(t)})_+-(m_{jj'\to h_{ii'jj'}}^{(t)})_+\right).
\end{align*}
and
\begin{multline*}
m_{ii'\to h_{ii'jj'}}^{(t+1)}=2\alpha w_{ii'}-\left\{\max_{k\neq i}(\alpha w_{ki'}+m_{ki'\to g_{i'}}^{(t)})\right\}_+-\left\{\max_{k'\neq i'}(\alpha w_{ik'}+m_{ik'\to f_{i}}^{(t)})\right\}_+\\
+\sum_{ii'kk'\neq ii'jj'}\left((\frac{\beta}{2}+m_{kk'\to h_{ii'kk'}}^{(t)})_+-(m_{kk'\to h_{ii'kk'}}^{(t)})_+\right).
\end{multline*}
We can simplify these equations further, by defining $m_{ii'\to f_i}^{(t)}\equiv w_{ii'}+m_{ii'\to f_i}^{(t)}$ and $m_{ii'\to g_{i'}}^{(t)}\equiv\alpha w_{ii'}+m_{ii'\to g_{i'}}^{(t)}$ and replacing $\tilde{\beta}$ with $\beta/2$ to obtain the following result.
\begin{lemma}
The max-product equations \eqref{eq:m(ij->fi)}-\eqref{eq:m(var->fun)} are equivalent to the simplified BP equations \eqref{eq:bp1}-\eqref{eq:bp2}.
\end{lemma}

}

\end{document}